\documentclass{amsart}
\usepackage[all]{xy}
\usepackage{tikz}
 \tikzstyle{int}=[circle, draw,fill=black,outer sep=0,minimum size=3pt, inner sep=0]
  \tikzstyle{ext}=[circle, draw=black,outer sep=0,inner sep=1pt]
\usepackage{latexsym}
\usepackage{amssymb}
\usepackage{amsmath}
\usepackage{amsthm}
\usepackage{amscd}
\usepackage{fancyhdr}
\usepackage{fullpage}
\usepackage{mathrsfs}
\input cyracc.def
\xyoption{arc}

\def\id{{\mbox{1 \hskip -8pt 1}}}
\newcommand{\sgn}{{\mathit s  \mathit g\mathit  n}}
 \newcommand{\lon}{\longrightarrow}
 \newcommand{\bu}{\bullet}
 
 \newcommand{\rar}{\rightarrow}

 \newcommand{\Z}{{\mathbb Z}}
 \newcommand{\bS}{{\mathbb S}}

 \newcommand{\R}{{\mathbb R}}
 
 \newcommand{\K}{{\mathbb K}}

\newcommand{\fGC}{\mathsf{fGC}}

 \newcommand{\ot}{\otimes}
 
 \newcommand{\Beq}{\begin{equation}}
 \newcommand{\Eeq}{\end{equation}}
 \newcommand{\Beqr}{\begin{eqnarray}}
 \newcommand{\Eeqr}{\end{eqnarray}}
 \newcommand{\Beqrn}{\begin{eqnarray*}}
 \newcommand{\Eeqrn}{\end{eqnarray*}}
 \newcommand{\Ba}{\begin{array}}
 \newcommand{\Ea}{\end{array}}
 \newcommand{\Bi}{\begin{itemize}}
 \newcommand{\Ei}{\end{itemize}}
 \newcommand{\Bc}{\begin{center}}
 \newcommand{\Ec}{\end{center}}

 \newcommand{\f}{{\mathcal O}}
 \newcommand{\cA}{{\mathcal A}}

 \newcommand{\cE}{{\mathcal E}}
 \newcommand{\cF}{{\mathcal F}}
 
 \newcommand{\caH}{{\mathcal H}}

 \newcommand{\caL}{{\mathcal L}}
 \newcommand{\cM}{{\mathcal M}}
 
 \newcommand{\cP}{{\mathcal P}}
 
 \newcommand{\cR}{{\mathcal R}}
 
 \newcommand{\cT}{{\mathcal T}}

 \newcommand{\cU}{{\mathcal U}}

 \newcommand{\Ga}{\Gamma}

 \newcommand{\Hom}{{\mathrm H\mathrm o\mathrm m}}
 
 \newcommand{\sip}{\smallskip}
 \newcommand{\bip}{\bigskip}
 \newcommand{\mip}{\vspace{2.5mm}}

 \newcommand{\fGCor}{\mathsf{fGC}^{or}}

 \newcommand{\LnB}{{\caL \mathit{ieb}_{odd}}}
  \newcommand{\LonB}{{\caL \mathit{ieb}_{odd}^\diamond}}
  \newcommand{\HoLonB}{{\caH \mathit{olieb}_{odd}^\diamond}}
  \newcommand{\HoLnB}{{\caH \mathit{olieb}_{odd}}}

  \newcommand{\LB}{\mathcal{L}\mathit{ieb}}
\newcommand{\LoB}{\mathcal{L}\mathit{ieb}^\diamond}

\newcommand{\HoLoB}{\mathcal{H}\mathit{olieb}^\diamond}

\newcommand{\wHoLonB}{\widehat{\mathcal{H}\mathit{olieb}}_{odd}^{\Ba{c}\vspace{-1mm}_{\hspace{-2mm}\diamond} \Ea}}

 \newcommand{\Der}{\mathrm{Der}}

 \newcommand{\gr}{\mathrm{gr}}
 
\theoremstyle{plain}
\swapnumbers

\newtheorem{prop-def}[theorem]{Proposition-definition}

\newtheorem{main-theorem}{Main~Theorem}[section]
\newtheorem{section-theorem}{Theorem}[section]
\newtheorem{section-corollary}{Corollary}[section]

\theoremstyle{definition}

\begin{document}

\sloppy

 \newenvironment{proo}{\begin{trivlist} \item{\sc {Proof.}}}
  {\hfill $\square$ \end{trivlist}}

\long\def\symbolfootnote[#1]#2{\begingroup%
\def\thefootnote{\fnsymbol{footnote}}\footnote[#1]{#2}\endgroup}

 \title{On quantizable  odd Lie bialgebras}

\author{Anton~Khoroshkin}
\address{Anton~Khoroshkin:
National Research University Higher School of Economics,
International Laboratory of Representation Theory and Mathematical Physics,
20 Myasnitskay Ulitsa, Moscow 101000, Russia
\newline
and
ITEP, Bolshaya Cheremushkinskaya 25, 117259, Moscow, Russia
}
\email{akhoroshkin@hse.ru}

\author{Sergei~Merkulov}
\address{Sergei~Merkulov:  Mathematics Research Unit, Luxembourg University,  Grand Duchy of Luxembourg}
\email{sergei.merkulov@uni.lu}

\author{Thomas~Willwacher}
\address{Thomas~Willwacher: Institute of Mathematics, University of Zurich, Zurich, Switzerland}
\email{thomas.willwacher@math.uzh.ch}

\begin{abstract} Motivated by the obstruction to the deformation quantization of Poisson structures in {\em infinite}\, dimensions we introduce the notion of a quantizable odd Lie bialgebra. The main result of the paper is a construction of the highly non-trivial
 minimal resolution of the properad governing such Lie bialgebras, and its link with the theory of so called {\em quantizable}\, Poisson structures.


\end{abstract}
 \maketitle

\section{Introduction}

 \subsection{Even and odd Lie bialgebras}
 A Lie ($c,d$)-bialgebra is a graded vector space $V$ which carries both a degree $c$ Lie algebra structure
$$
\begin{xy}
 <0mm,0.66mm>*{};<0mm,3mm>*{}**@{-},
 <0.39mm,-0.39mm>*{};<2.2mm,-2.2mm>*{}**@{-},
 <-0.35mm,-0.35mm>*{};<-2.2mm,-2.2mm>*{}**@{-},
 <0mm,0mm>*{\circ};<0mm,0mm>*{}**@{},
   <0mm,0.66mm>*{};<0mm,3.4mm>*{^1}**@{},
   <0.39mm,-0.39mm>*{};<2.9mm,-4mm>*{^2}**@{},
   <-0.35mm,-0.35mm>*{};<-2.8mm,-4mm>*{^1}**@{},
\end{xy}\simeq  [\ ,\ ]: V\ot V\rar V[c]
$$
and a degree $d$ Lie coalgebra structure,
$$
\begin{xy}
 <0mm,-0.55mm>*{};<0mm,-2.5mm>*{}**@{-},
 <0.5mm,0.5mm>*{};<2.2mm,2.2mm>*{}**@{-},
 <-0.48mm,0.48mm>*{};<-2.2mm,2.2mm>*{}**@{-},
 <0mm,0mm>*{\circ};<0mm,0mm>*{}**@{},
 <0mm,-0.55mm>*{};<0mm,-3.8mm>*{_1}**@{},
 <0.5mm,0.5mm>*{};<2.7mm,2.8mm>*{^2}**@{},
 <-0.48mm,0.48mm>*{};<-2.7mm,2.8mm>*{^1}**@{},
 \end{xy}\simeq \Delta:V \rar V\ot V[d]
$$
satisfying the following compatibility condition:
\begin{multline*}
\begin{xy}
<0mm,2.47mm>*{};<0mm,0.12mm>*{}**@{-},
<0.5mm,3.5mm>*{};<2.2mm,5.2mm>*{}**@{-},
<-0.48mm,3.48mm>*{};<-2.2mm,5.2mm>*{}**@{-},
<0mm,3mm>*{\circ};<0mm,3mm>*{}**@{},
<0mm,-0.8mm>*{\circ};<0mm,-0.8mm>*{}**@{},
<-0.39mm,-1.2mm>*{};<-2.2mm,-3.5mm>*{}**@{-},
<0.39mm,-1.2mm>*{};<2.2mm,-3.5mm>*{}**@{-},
<0.5mm,3.5mm>*{};<2.8mm,5.7mm>*{^2}**@{},
<-0.48mm,3.48mm>*{};<-2.8mm,5.7mm>*{^1}**@{},
<0mm,-0.8mm>*{};<-2.7mm,-5.2mm>*{^1}**@{},
<0mm,-0.8mm>*{};<2.7mm,-5.2mm>*{^2}**@{},
\end{xy}
\  - \
\begin{xy}
<0mm,-1.3mm>*{};<0mm,-3.5mm>*{}**@{-},
<0.38mm,-0.2mm>*{};<2.0mm,2.0mm>*{}**@{-},
<-0.38mm,-0.2mm>*{};<-2.2mm,2.2mm>*{}**@{-},
<0mm,-0.8mm>*{\circ};<0mm,0.8mm>*{}**@{},
<2.4mm,2.4mm>*{\circ};<2.4mm,2.4mm>*{}**@{},
<2.77mm,2.0mm>*{};<4.4mm,-0.8mm>*{}**@{-},
<2.4mm,3mm>*{};<2.4mm,5.2mm>*{}**@{-},
<0mm,-1.3mm>*{};<0mm,-5.3mm>*{^1}**@{},
<2.5mm,2.3mm>*{};<5.1mm,-2.6mm>*{^2}**@{},
<2.4mm,2.5mm>*{};<2.4mm,5.7mm>*{^2}**@{},
<-0.38mm,-0.2mm>*{};<-2.8mm,2.5mm>*{^1}**@{},
\end{xy}
\  - (-1)^{d}\
\begin{xy}
<0mm,-1.3mm>*{};<0mm,-3.5mm>*{}**@{-},
<0.38mm,-0.2mm>*{};<2.0mm,2.0mm>*{}**@{-},
<-0.38mm,-0.2mm>*{};<-2.2mm,2.2mm>*{}**@{-},
<0mm,-0.8mm>*{\circ};<0mm,0.8mm>*{}**@{},
<2.4mm,2.4mm>*{\circ};<2.4mm,2.4mm>*{}**@{},
<2.77mm,2.0mm>*{};<4.4mm,-0.8mm>*{}**@{-},
<2.4mm,3mm>*{};<2.4mm,5.2mm>*{}**@{-},
<0mm,-1.3mm>*{};<0mm,-5.3mm>*{^2}**@{},
<2.5mm,2.3mm>*{};<5.1mm,-2.6mm>*{^1}**@{},
<2.4mm,2.5mm>*{};<2.4mm,5.7mm>*{^2}**@{},
<-0.38mm,-0.2mm>*{};<-2.8mm,2.5mm>*{^1}**@{},
\end{xy}
\  -(-1)^{c+d} \
\begin{xy}
<0mm,-1.3mm>*{};<0mm,-3.5mm>*{}**@{-},
<0.38mm,-0.2mm>*{};<2.0mm,2.0mm>*{}**@{-},
<-0.38mm,-0.2mm>*{};<-2.2mm,2.2mm>*{}**@{-},
<0mm,-0.8mm>*{\circ};<0mm,0.8mm>*{}**@{},
<2.4mm,2.4mm>*{\circ};<2.4mm,2.4mm>*{}**@{},
<2.77mm,2.0mm>*{};<4.4mm,-0.8mm>*{}**@{-},
<2.4mm,3mm>*{};<2.4mm,5.2mm>*{}**@{-},
<0mm,-1.3mm>*{};<0mm,-5.3mm>*{^2}**@{},
<2.5mm,2.3mm>*{};<5.1mm,-2.6mm>*{^1}**@{},
<2.4mm,2.5mm>*{};<2.4mm,5.7mm>*{^1}**@{},
<-0.38mm,-0.2mm>*{};<-2.8mm,2.5mm>*{^2}**@{},
\end{xy}
\ - (-1)^{c} \
\begin{xy}
<0mm,-1.3mm>*{};<0mm,-3.5mm>*{}**@{-},
<0.38mm,-0.2mm>*{};<2.0mm,2.0mm>*{}**@{-},
<-0.38mm,-0.2mm>*{};<-2.2mm,2.2mm>*{}**@{-},
<0mm,-0.8mm>*{\circ};<0mm,0.8mm>*{}**@{},
<2.4mm,2.4mm>*{\circ};<2.4mm,2.4mm>*{}**@{},
<2.77mm,2.0mm>*{};<4.4mm,-0.8mm>*{}**@{-},
<2.4mm,3mm>*{};<2.4mm,5.2mm>*{}**@{-},
<0mm,-1.3mm>*{};<0mm,-5.3mm>*{^1}**@{},
<2.5mm,2.3mm>*{};<5.1mm,-2.6mm>*{^2}**@{},
<2.4mm,2.5mm>*{};<2.4mm,5.7mm>*{^1}**@{},
<-0.38mm,-0.2mm>*{};<-2.8mm,2.5mm>*{^2}**@{},
\end{xy}\ \simeq \\
\simeq\
\Delta([x_1,x_2]) - [\Delta(x_1),x_2\otimes 1 + 1\otimes x_2] \pm [x_1\otimes 1 + 1\otimes x_1,\Delta(x_2)] = 0.
\end{multline*}
 If the $\Z_2$-parities
of both structures are the same, i.e.\ if $c+d\in 2\Z$, the Lie bialgebra is called {\em even}, if the $\Z_2$-parities
are opposite, $c+d\in 2\Z+1$, it is called {\em odd}.

\sip

In the even case the most interesting for applications Lie bialgebras have $c=d=0$. Such Lie bialgebras were introduced by Drinfeld in \cite{D1}  in the context of the
theory of Yang-Baxter equations, and they have since found numerous applications, most prominently in the theory
of Hopf algebra deformations of universal enveloping algebras (see the book \cite{ES} and references cited there).
If the composition of the cobracket and bracket of a Lie bialgebra is zero, that is
\Beq\label{1: involutivity condition}
\Ba{c}\resizebox{4mm}{!}
{\xy
 (0,0)*{\circ}="a",
(0,6)*{\circ}="b",
(3,3)*{}="c",
(-3,3)*{}="d",
 (0,9)*{}="b'",
(0,-3)*{}="a'",
\ar@{-} "a";"c" <0pt>
\ar @{-} "a";"d" <0pt>
\ar @{-} "a";"a'" <0pt>
\ar @{-} "b";"c" <0pt>
\ar @{-} "b";"d" <0pt>
\ar @{-} "b";"b'" <0pt>
\endxy}
\Ea = 0,
\Eeq
then the Lie bialgebra is called \emph{involutive}.
This additional constraint is satisfied in many interesting examples studied in homological algebra, string topology, symplectic field theory,  Lagrangian Floer theory of higher genus, and the theory of cohomology groups $H(\cM_{g,n})$ of moduli spaces of algebraic curves with labelings of punctures skewsymmetrized \cite{D1,ES,Ch,ChSu,CFL,S,CMW,MW0}.
\sip

In the odd case the most interesting for applications Lie bialgebras have $c=1$, $d=0$.
They have been introduced in \cite{Me1} and have seen applications in Poisson geometry, deformation quantization
of Poisson structures \cite{Me2} and in the theory of cohomology groups $H(\cM_{g,n})$ of moduli spaces of algebraic curves with labelings of punctures { symmetrized} \cite{MW0}.

\sip

The homotopy and deformation theories of even/odd Lie bialgebras and also of involutive Lie bialgebras have been studied in \cite{CMW,MW}. A key tool in those studies is a minimal resolution of the properad governing the algebraic structure under consideration.

\sip

The minimal resolutions of properads $\LB$ and $\LnB$ governing
even and, respectively, odd  Lie bialgebras were constructed in \cite{Ko,MaVo,V} and, respectively, in \cite{Me1,Me2}. Constructing a minimal resolution
$\HoLoB$ of the properad  $\LoB$ governing {\em involutive}\, Lie bialgebras turned out to be a more difficult problem, and that goal was achieved only very recently in \cite{CMW}.

\subsection{Quantizable odd  Lie bialgebras} For odd Lie bialgebras the involutivity condition
(\ref{1: involutivity condition}) is trivial, i.e.\ it is satisfied automatically for any odd Lie bialgebra $V$. There is, however, a higher genus analogue of  that condition,
\Beq\label{1: quantizability condition}
\Ba{c}
\resizebox{8mm}{!}
{
\begin{xy}
 <0mm,-1.3mm>*{};<0mm,-3.5mm>*{}**@{-},
 <0.4mm,0.0mm>*{};<2.4mm,2.1mm>*{}**@{-},
 <-0.38mm,-0.2mm>*{};<-2.8mm,2.5mm>*{}**@{-},
<0mm,0.8mm>*{}**@{},
<2.45mm,2.35mm>*{}**@{},
 <2.4mm,2.8mm>*{};<0mm,5mm>*{}**@{-},
 <3.35mm,2.9mm>*{};<5.5mm,5mm>*{}**@{-},
<-2.8mm,2.5mm>*{};<0mm,5mm>*{}**@{-},
<2.96mm,5mm>*{};
<0.2mm,5.1mm>*{};<2.8mm,7.5mm>*{}**@{-},
<0.2mm,5.1mm>*{};<2.8mm,7.5mm>*{}**@{-},
<5.5mm,5mm>*{};<3.4mm,7.0mm>*{}**@{-},
<2.9mm,8.2mm>*{};<2.9mm,10.5mm>*{}**@{-},
<-2.8mm,2.5mm>*{};<0mm,5mm>*{\circ}**@{},
<0mm,-0.8mm>*{\circ};
 <2.96mm,2.4mm>*{\circ};
 <2.96mm,7.5mm>*{\circ}**@{},
    \end{xy}}
\Ea\ =\ 0,
\Eeq
which is highly non-trivial, and which can be considered as an odd analogue of (\ref{1: involutivity condition}). We prefer, however, to call odd Lie bialgebras satisfying the extra constraint (\ref{1: quantizability condition}) {\em quantizable}\, ones rather than involutive.
The reason for this terminology is explained in \S 2: this graph controls, on the one hand, the obstruction to the universal quantization of Poisson structures in {\em infinite}\, dimensions \cite{Me1,Me2}, and, on the other hand, it controls the obstruction to the existence of  geometrically
meaningful Frobenius structure on $\mathrm{Chains}(\R)$ \cite{JF1}.

\sip

Our main result is an explicit construction in \S 3 of a (highly non-obvious) minimal resolution $\HoLonB$ of the properad $\LonB$
governing quantizable Lie bialgebras. One of the key tricks in \cite{CMW} used to solve an analogous problem for the properad $\LoB$ of even involutive Lie bialgebras  reduced the ``hard" problem of computing
the cohomology of some dg properad to the  ``easy" computation of the minimal resolutions of a family of some auxiliary {\em quadratic}\, algebras. Remarkable enough,  this approach works for the  constraint (\ref{1: quantizability condition}) as well,
but it leads instead to a certain  family of {\em cubic}\, homogeneous algebras which are studied in the Appendix.
\sip

Another important technical ingredient in our construction of $\HoLonB$ comes from the
paper \cite{MW}, in which the cohomologies of the deformation complexes of the properads
 $\LB$, $\LnB$ and $\LoB$ have been computed, and it was proven in particular that the properad $\LnB$
admits precisely one non-trivial deformation; in fact it is that unique non-trivial deformation which leads us to the dg properad  $\HoLonB$. We explain this link in \S 3.

\sip

There are plenty of examples of $\HoLonB$ algebra structures associated with ordinary (or formal power series) Poisson structures $\pi\in \cT_{poly}(\R^d)$ on $\R^d$, $d\geq 2$, which vanish at $0\in \R^d$; for a generic Poisson structure $\pi$ on $\R^d$ the associated $\HoLonB$ algebra structure $\pi^\diamond$,
$$
\pi^\diamond: \HoLonB \lon \cE nd_{\R^d},
$$
in $\R^d$
is highly non-trivial and can be given explicitly only by transcendental formulae (i.e.\ the ones involving converging integrals over suitable configuration spaces \cite{MW2}). Such a $\HoLonB$ algebra structure
$\pi^\diamond$ in $\R^d$ can also be interpreted as a formal power series bivector field $\pi^\diamond\in \cT_{poly}(\R^d)[[\hbar]]$,
$$
\pi^\diamond = \pi_0 + \hbar \pi_1 + \hbar^3 \pi_2 +\ldots,  \ \ \ \ \ \
$$
satisfying a certain formal power series equation,
$$
\frac{1}{2}[\pi^\diamond,\pi^\diamond]_2 + \frac{\hbar}{4!} [\pi^\diamond,\pi^\diamond,\pi^\diamond,\pi^\diamond]_4
+ \ldots =0,
$$
Here the collection of operators,
$$
\left\{[\, \ , \ldots ,\ ]_{2n}: \cT_{poly}(\R^d)^{\ot 2n}\rar \cT_{poly}(\R^d)[3-4n] \right\}_{n\geq 1}
$$
defines a so called {\em Kontsevich-Shoikhet} $\caL ie_\infty$ structure \cite{Sh} in  $\cT_{poly}(\R^d)$
with $[\ ,\ ]_2$ being the standard Schouten bracket.
This huge class $\{\pi^\diamond\}$ of highly-non-trivial examples of $\HoLonB$ algebra structures in $\R^d$ depends in general on the choice of Drinfeld associator \cite{MW2}; it motivated much our current study
of the homotopy theory of odd Lie bialgebras.

\subsection{Some notation}
 The set $\{1,2, \ldots, n\}$ is abbreviated to $[n]$;  its group of automorphisms is
denoted by $\bS_n$;
the trivial one-dimensional representation of
 $\bS_n$ is denoted by $\id_n$, while its one dimensional sign representation is
 denoted by $\sgn_n$.
The cardinality of a finite set $A$ is denoted by $\# A$.

\sip

We work throughout in the category of $\Z$-graded vector spaces over a field $\K$
of characteristic zero.
If $V=\oplus_{i\in \Z} V^i$ is a graded vector space, then
$V[k]$ stands for the graded vector space with $V[k]^i:=V^{i+k}$. For a
prop(erad) $\cP$ we denote by $\cP\{k\}$ a prop(erad) which is uniquely defined by
 the following property:
for any graded vector space $V$ a representation
of $\cP\{k\}$ in $V$ is identical to a representation of  $\cP$ in $V[k]$.


\subsection*{Acknowledgements} We are very grateful to Jan-Erik Roos for valuable correspondences concerning Koszulness theory of homogeneous algebras and to Dmitri Piontkovski for inspiring discussions on algebras of homological dimension $2$ and strongly free Gr\"obner bases.
\sip

A.Kh. acknowledges partial support by RFBR grant 15-01-09242,  by Dynasty foundation and Simons-IUM fellowship, partially
 prepared within the framework of a subsidy granted to the National Research University Higher School of Economics by the Government of the Russian Federation for the implementation of the Global Competitiveness Program.
S.M. acknowledges support
by the Swedish Vetenskapr\aa det, grant 2012-5478.
T.W. acknowledges partial support by the Swiss National Science Foundation, grant 200021\_150012, and by the NCCR SwissMAP of the Swiss National Science Foundation.

\bip

{
\Large
 \section{\bf Quantizable odd Lie bialgebras}
 }

 \sip

 \subsection{Odd lie bialgebras}

By definition \cite{Me1}, the properad, $\LonB$, of odd  Lie bialgebras is a quadratic properad
given as the quotient,
$$
\LonB:=\cF ree\langle E\rangle/\langle\cR \rangle,
$$
of the free properad generated by an  $\bS$-bimodule $E=\{E(m,n)\}_{m,n\geq 1}$ with
 all $E(m,n)=0$ except
$$
E(2,1):= \sgn_2\ot \id_1=\mbox{span}\left\langle
\begin{xy}
 <0mm,-0.55mm>*{};<0mm,-2.5mm>*{}**@{-},
 <0.5mm,0.5mm>*{};<2.2mm,2.2mm>*{}**@{-},
 <-0.48mm,0.48mm>*{};<-2.2mm,2.2mm>*{}**@{-},
 <0mm,0mm>*{\circ};<0mm,0mm>*{}**@{},
 <0mm,-0.55mm>*{};<0mm,-3.8mm>*{_1}**@{},
 <0.5mm,0.5mm>*{};<2.7mm,2.8mm>*{^2}**@{},
 <-0.48mm,0.48mm>*{};<-2.7mm,2.8mm>*{^1}**@{},
 \end{xy}
=-
\begin{xy}
 <0mm,-0.55mm>*{};<0mm,-2.5mm>*{}**@{-},
 <0.5mm,0.5mm>*{};<2.2mm,2.2mm>*{}**@{-},
 <-0.48mm,0.48mm>*{};<-2.2mm,2.2mm>*{}**@{-},
 <0mm,0mm>*{\circ};<0mm,0mm>*{}**@{},
 <0mm,-0.55mm>*{};<0mm,-3.8mm>*{_1}**@{},
 <0.5mm,0.5mm>*{};<2.7mm,2.8mm>*{^1}**@{},
 <-0.48mm,0.48mm>*{};<-2.7mm,2.8mm>*{^2}**@{},
 \end{xy}
   \right\rangle
$$
$$
E(1,2):= \id_1\ot \id_2[-1]=\mbox{span}\left\langle
\begin{xy}
 <0mm,0.66mm>*{};<0mm,3mm>*{}**@{-},
 <0.39mm,-0.39mm>*{};<2.2mm,-2.2mm>*{}**@{-},
 <-0.35mm,-0.35mm>*{};<-2.2mm,-2.2mm>*{}**@{-},
 <0mm,0mm>*{\circ};<0mm,0mm>*{}**@{},
   <0mm,0.66mm>*{};<0mm,3.4mm>*{^1}**@{},
   <0.39mm,-0.39mm>*{};<2.9mm,-4mm>*{^2}**@{},
   <-0.35mm,-0.35mm>*{};<-2.8mm,-4mm>*{^1}**@{},
\end{xy}=
\begin{xy}
 <0mm,0.66mm>*{};<0mm,3mm>*{}**@{-},
 <0.39mm,-0.39mm>*{};<2.2mm,-2.2mm>*{}**@{-},
 <-0.35mm,-0.35mm>*{};<-2.2mm,-2.2mm>*{}**@{-},
 <0mm,0mm>*{\circ};<0mm,0mm>*{}**@{},
   <0mm,0.66mm>*{};<0mm,3.4mm>*{^1}**@{},
   <0.39mm,-0.39mm>*{};<2.9mm,-4mm>*{^1}**@{},
   <-0.35mm,-0.35mm>*{};<-2.8mm,-4mm>*{^2}**@{},
\end{xy}
\right\rangle
$$
modulo the ideal generated by the following relations
\Beq\label{R for LieB}
\cR:\left\{
\Ba{c}
\begin{xy}
 <0mm,0mm>*{\circ};<0mm,0mm>*{}**@{},
 <0mm,-0.49mm>*{};<0mm,-3.0mm>*{}**@{-},
 <0.49mm,0.49mm>*{};<1.9mm,1.9mm>*{}**@{-},
 <-0.5mm,0.5mm>*{};<-1.9mm,1.9mm>*{}**@{-},
 <-2.3mm,2.3mm>*{\circ};<-2.3mm,2.3mm>*{}**@{},
 <-1.8mm,2.8mm>*{};<0mm,4.9mm>*{}**@{-},
 <-2.8mm,2.9mm>*{};<-4.6mm,4.9mm>*{}**@{-},
   <0.49mm,0.49mm>*{};<2.7mm,2.3mm>*{^3}**@{},
   <-1.8mm,2.8mm>*{};<0.4mm,5.3mm>*{^2}**@{},
   <-2.8mm,2.9mm>*{};<-5.1mm,5.3mm>*{^1}**@{},
 \end{xy}
\ + \
\begin{xy}
 <0mm,0mm>*{\circ};<0mm,0mm>*{}**@{},
 <0mm,-0.49mm>*{};<0mm,-3.0mm>*{}**@{-},
 <0.49mm,0.49mm>*{};<1.9mm,1.9mm>*{}**@{-},
 <-0.5mm,0.5mm>*{};<-1.9mm,1.9mm>*{}**@{-},
 <-2.3mm,2.3mm>*{\circ};<-2.3mm,2.3mm>*{}**@{},
 <-1.8mm,2.8mm>*{};<0mm,4.9mm>*{}**@{-},
 <-2.8mm,2.9mm>*{};<-4.6mm,4.9mm>*{}**@{-},
   <0.49mm,0.49mm>*{};<2.7mm,2.3mm>*{^2}**@{},
   <-1.8mm,2.8mm>*{};<0.4mm,5.3mm>*{^1}**@{},
   <-2.8mm,2.9mm>*{};<-5.1mm,5.3mm>*{^3}**@{},
 \end{xy}
\ + \
\begin{xy}
 <0mm,0mm>*{\circ};<0mm,0mm>*{}**@{},
 <0mm,-0.49mm>*{};<0mm,-3.0mm>*{}**@{-},
 <0.49mm,0.49mm>*{};<1.9mm,1.9mm>*{}**@{-},
 <-0.5mm,0.5mm>*{};<-1.9mm,1.9mm>*{}**@{-},
 <-2.3mm,2.3mm>*{\circ};<-2.3mm,2.3mm>*{}**@{},
 <-1.8mm,2.8mm>*{};<0mm,4.9mm>*{}**@{-},
 <-2.8mm,2.9mm>*{};<-4.6mm,4.9mm>*{}**@{-},
   <0.49mm,0.49mm>*{};<2.7mm,2.3mm>*{^1}**@{},
   <-1.8mm,2.8mm>*{};<0.4mm,5.3mm>*{^3}**@{},
   <-2.8mm,2.9mm>*{};<-5.1mm,5.3mm>*{^2}**@{},
 \end{xy}\ =\ 0,
 \vspace{3mm}\\
 \begin{xy}
 <0mm,0mm>*{\circ};<0mm,0mm>*{}**@{},
 <0mm,0.69mm>*{};<0mm,3.0mm>*{}**@{-},
 <0.39mm,-0.39mm>*{};<2.4mm,-2.4mm>*{}**@{-},
 <-0.35mm,-0.35mm>*{};<-1.9mm,-1.9mm>*{}**@{-},
 <-2.4mm,-2.4mm>*{\circ};<-2.4mm,-2.4mm>*{}**@{},
 <-2.0mm,-2.8mm>*{};<0mm,-4.9mm>*{}**@{-},
 <-2.8mm,-2.9mm>*{};<-4.7mm,-4.9mm>*{}**@{-},
    <0.39mm,-0.39mm>*{};<3.3mm,-4.0mm>*{^3}**@{},
    <-2.0mm,-2.8mm>*{};<0.5mm,-6.7mm>*{^2}**@{},
    <-2.8mm,-2.9mm>*{};<-5.2mm,-6.7mm>*{^1}**@{},
 \end{xy}
\ + \
 \begin{xy}
 <0mm,0mm>*{\circ};<0mm,0mm>*{}**@{},
 <0mm,0.69mm>*{};<0mm,3.0mm>*{}**@{-},
 <0.39mm,-0.39mm>*{};<2.4mm,-2.4mm>*{}**@{-},
 <-0.35mm,-0.35mm>*{};<-1.9mm,-1.9mm>*{}**@{-},
 <-2.4mm,-2.4mm>*{\circ};<-2.4mm,-2.4mm>*{}**@{},
 <-2.0mm,-2.8mm>*{};<0mm,-4.9mm>*{}**@{-},
 <-2.8mm,-2.9mm>*{};<-4.7mm,-4.9mm>*{}**@{-},
    <0.39mm,-0.39mm>*{};<3.3mm,-4.0mm>*{^2}**@{},
    <-2.0mm,-2.8mm>*{};<0.5mm,-6.7mm>*{^1}**@{},
    <-2.8mm,-2.9mm>*{};<-5.2mm,-6.7mm>*{^3}**@{},
 \end{xy}
\ + \
 \begin{xy}
 <0mm,0mm>*{\circ};<0mm,0mm>*{}**@{},
 <0mm,0.69mm>*{};<0mm,3.0mm>*{}**@{-},
 <0.39mm,-0.39mm>*{};<2.4mm,-2.4mm>*{}**@{-},
 <-0.35mm,-0.35mm>*{};<-1.9mm,-1.9mm>*{}**@{-},
 <-2.4mm,-2.4mm>*{\circ};<-2.4mm,-2.4mm>*{}**@{},
 <-2.0mm,-2.8mm>*{};<0mm,-4.9mm>*{}**@{-},
 <-2.8mm,-2.9mm>*{};<-4.7mm,-4.9mm>*{}**@{-},
    <0.39mm,-0.39mm>*{};<3.3mm,-4.0mm>*{^1}**@{},
    <-2.0mm,-2.8mm>*{};<0.5mm,-6.7mm>*{^3}**@{},
    <-2.8mm,-2.9mm>*{};<-5.2mm,-6.7mm>*{^2}**@{},
 \end{xy}\ =\ 0,
 \\
 \begin{xy}
 <0mm,2.47mm>*{};<0mm,0.12mm>*{}**@{-},
 <0.5mm,3.5mm>*{};<2.2mm,5.2mm>*{}**@{-},
 <-0.48mm,3.48mm>*{};<-2.2mm,5.2mm>*{}**@{-},
 <0mm,3mm>*{\circ};<0mm,3mm>*{}**@{},
  <0mm,-0.8mm>*{\circ};<0mm,-0.8mm>*{}**@{},
<-0.39mm,-1.2mm>*{};<-2.2mm,-3.5mm>*{}**@{-},
 <0.39mm,-1.2mm>*{};<2.2mm,-3.5mm>*{}**@{-},
     <0.5mm,3.5mm>*{};<2.8mm,5.7mm>*{^2}**@{},
     <-0.48mm,3.48mm>*{};<-2.8mm,5.7mm>*{^1}**@{},
   <0mm,-0.8mm>*{};<-2.7mm,-5.2mm>*{^1}**@{},
   <0mm,-0.8mm>*{};<2.7mm,-5.2mm>*{^2}**@{},
\end{xy}
\  - \
\begin{xy}
 <0mm,-1.3mm>*{};<0mm,-3.5mm>*{}**@{-},
 <0.38mm,-0.2mm>*{};<2.0mm,2.0mm>*{}**@{-},
 <-0.38mm,-0.2mm>*{};<-2.2mm,2.2mm>*{}**@{-},
<0mm,-0.8mm>*{\circ};<0mm,0.8mm>*{}**@{},
 <2.4mm,2.4mm>*{\circ};<2.4mm,2.4mm>*{}**@{},
 <2.77mm,2.0mm>*{};<4.4mm,-0.8mm>*{}**@{-},
 <2.4mm,3mm>*{};<2.4mm,5.2mm>*{}**@{-},
     <0mm,-1.3mm>*{};<0mm,-5.3mm>*{^1}**@{},
     <2.5mm,2.3mm>*{};<5.1mm,-2.6mm>*{^2}**@{},
    <2.4mm,2.5mm>*{};<2.4mm,5.7mm>*{^2}**@{},
    <-0.38mm,-0.2mm>*{};<-2.8mm,2.5mm>*{^1}**@{},
    \end{xy}
\  -\
\begin{xy}
 <0mm,-1.3mm>*{};<0mm,-3.5mm>*{}**@{-},
 <0.38mm,-0.2mm>*{};<2.0mm,2.0mm>*{}**@{-},
 <-0.38mm,-0.2mm>*{};<-2.2mm,2.2mm>*{}**@{-},
<0mm,-0.8mm>*{\circ};<0mm,0.8mm>*{}**@{},
 <2.4mm,2.4mm>*{\circ};<2.4mm,2.4mm>*{}**@{},
 <2.77mm,2.0mm>*{};<4.4mm,-0.8mm>*{}**@{-},
 <2.4mm,3mm>*{};<2.4mm,5.2mm>*{}**@{-},
     <0mm,-1.3mm>*{};<0mm,-5.3mm>*{^2}**@{},
     <2.5mm,2.3mm>*{};<5.1mm,-2.6mm>*{^1}**@{},
    <2.4mm,2.5mm>*{};<2.4mm,5.7mm>*{^2}**@{},
    <-0.38mm,-0.2mm>*{};<-2.8mm,2.5mm>*{^1}**@{},
    \end{xy}
\  + \
\begin{xy}
 <0mm,-1.3mm>*{};<0mm,-3.5mm>*{}**@{-},
 <0.38mm,-0.2mm>*{};<2.0mm,2.0mm>*{}**@{-},
 <-0.38mm,-0.2mm>*{};<-2.2mm,2.2mm>*{}**@{-},
<0mm,-0.8mm>*{\circ};<0mm,0.8mm>*{}**@{},
 <2.4mm,2.4mm>*{\circ};<2.4mm,2.4mm>*{}**@{},
 <2.77mm,2.0mm>*{};<4.4mm,-0.8mm>*{}**@{-},
 <2.4mm,3mm>*{};<2.4mm,5.2mm>*{}**@{-},
     <0mm,-1.3mm>*{};<0mm,-5.3mm>*{^2}**@{},
     <2.5mm,2.3mm>*{};<5.1mm,-2.6mm>*{^1}**@{},
    <2.4mm,2.5mm>*{};<2.4mm,5.7mm>*{^1}**@{},
    <-0.38mm,-0.2mm>*{};<-2.8mm,2.5mm>*{^2}**@{},
    \end{xy}
\ + \
\begin{xy}
 <0mm,-1.3mm>*{};<0mm,-3.5mm>*{}**@{-},
 <0.38mm,-0.2mm>*{};<2.0mm,2.0mm>*{}**@{-},
 <-0.38mm,-0.2mm>*{};<-2.2mm,2.2mm>*{}**@{-},
<0mm,-0.8mm>*{\circ};<0mm,0.8mm>*{}**@{},
 <2.4mm,2.4mm>*{\circ};<2.4mm,2.4mm>*{}**@{},
 <2.77mm,2.0mm>*{};<4.4mm,-0.8mm>*{}**@{-},
 <2.4mm,3mm>*{};<2.4mm,5.2mm>*{}**@{-},
     <0mm,-1.3mm>*{};<0mm,-5.3mm>*{^1}**@{},
     <2.5mm,2.3mm>*{};<5.1mm,-2.6mm>*{^2}**@{},
    <2.4mm,2.5mm>*{};<2.4mm,5.7mm>*{^1}**@{},
    <-0.38mm,-0.2mm>*{};<-2.8mm,2.5mm>*{^2}**@{},
    \end{xy}\ =\ 0.
\Ea
\right.
\Eeq
A minimal resolution $\HoLnB$ of $\LnB$ was constructed in \cite{Me1,Me2}. It is a free properad,
$$
\HoLnB=\cF ree \langle\hat{E}\rangle
$$
generated by an $\bS$--bimodule $\hat{E}=\{ \hat{E}(m,n)\}_{m,n\geq 1, m+n\geq 3}$,
$$
\hat{E}(m,n):= sgn_m\ot \id_n[m-2]=\mbox{span}\left\langle
\resizebox{14mm}{!}{\begin{xy}
 <0mm,0mm>*{\circ};<0mm,0mm>*{}**@{},
 <-0.6mm,0.44mm>*{};<-8mm,5mm>*{}**@{-},
 <-0.4mm,0.7mm>*{};<-4.5mm,5mm>*{}**@{-},
 <0mm,0mm>*{};<-1mm,5mm>*{\ldots}**@{},
 <0.4mm,0.7mm>*{};<4.5mm,5mm>*{}**@{-},
 <0.6mm,0.44mm>*{};<8mm,5mm>*{}**@{-},
   <0mm,0mm>*{};<-8.5mm,5.5mm>*{^1}**@{},
   <0mm,0mm>*{};<-5mm,5.5mm>*{^2}**@{},
   <0mm,0mm>*{};<4.5mm,5.5mm>*{^{m\hspace{-0.5mm}-\hspace{-0.5mm}1}}**@{},
   <0mm,0mm>*{};<9.0mm,5.5mm>*{^m}**@{},
 <-0.6mm,-0.44mm>*{};<-8mm,-5mm>*{}**@{-},
 <-0.4mm,-0.7mm>*{};<-4.5mm,-5mm>*{}**@{-},
 <0mm,0mm>*{};<-1mm,-5mm>*{\ldots}**@{},
 <0.4mm,-0.7mm>*{};<4.5mm,-5mm>*{}**@{-},
 <0.6mm,-0.44mm>*{};<8mm,-5mm>*{}**@{-},
   <0mm,0mm>*{};<-8.5mm,-6.9mm>*{^1}**@{},
   <0mm,0mm>*{};<-5mm,-6.9mm>*{^2}**@{},
   <0mm,0mm>*{};<4.5mm,-6.9mm>*{^{n\hspace{-0.5mm}-\hspace{-0.5mm}1}}**@{},
   <0mm,0mm>*{};<9.0mm,-6.9mm>*{^n}**@{},
 \end{xy}}
\right\rangle,
$$
and comes equipped with the differential
$$
\delta
\resizebox{14mm}{!}{\begin{xy}
 <0mm,0mm>*{\circ};<0mm,0mm>*{}**@{},
 <-0.6mm,0.44mm>*{};<-8mm,5mm>*{}**@{-},
 <-0.4mm,0.7mm>*{};<-4.5mm,5mm>*{}**@{-},
 <0mm,0mm>*{};<-1mm,5mm>*{\ldots}**@{},
 <0.4mm,0.7mm>*{};<4.5mm,5mm>*{}**@{-},
 <0.6mm,0.44mm>*{};<8mm,5mm>*{}**@{-},
   <0mm,0mm>*{};<-8.5mm,5.5mm>*{^1}**@{},
   <0mm,0mm>*{};<-5mm,5.5mm>*{^2}**@{},
   <0mm,0mm>*{};<4.5mm,5.5mm>*{^{m\hspace{-0.5mm}-\hspace{-0.5mm}1}}**@{},
   <0mm,0mm>*{};<9.0mm,5.5mm>*{^m}**@{},
 <-0.6mm,-0.44mm>*{};<-8mm,-5mm>*{}**@{-},
 <-0.4mm,-0.7mm>*{};<-4.5mm,-5mm>*{}**@{-},
 <0mm,0mm>*{};<-1mm,-5mm>*{\ldots}**@{},
 <0.4mm,-0.7mm>*{};<4.5mm,-5mm>*{}**@{-},
 <0.6mm,-0.44mm>*{};<8mm,-5mm>*{}**@{-},
   <0mm,0mm>*{};<-8.5mm,-6.9mm>*{^1}**@{},
   <0mm,0mm>*{};<-5mm,-6.9mm>*{^2}**@{},
   <0mm,0mm>*{};<4.5mm,-6.9mm>*{^{n\hspace{-0.5mm}-\hspace{-0.5mm}1}}**@{},
   <0mm,0mm>*{};<9.0mm,-6.9mm>*{^n}**@{},
 \end{xy}}
\ \ = \ \
 \sum_{[1,\ldots,m]=I_1\sqcup I_2\atop
 {|I_1|\geq 0, |I_2|\geq 1}}
 \sum_{[1,\ldots,n]=J_1\sqcup J_2\atop
 {|J_1|\geq 1, |J_2|\geq 1}
}\hspace{0mm}
\pm\ \resizebox{24mm}{!}{ \begin{xy}
 <0mm,0mm>*{\circ};<0mm,0mm>*{}**@{},
 <-0.6mm,0.44mm>*{};<-8mm,5mm>*{}**@{-},
 <-0.4mm,0.7mm>*{};<-4.5mm,5mm>*{}**@{-},
 <0mm,0mm>*{};<0mm,5mm>*{\ldots}**@{},
 <0.4mm,0.7mm>*{};<4.5mm,5mm>*{}**@{-},
 <0.6mm,0.44mm>*{};<12.4mm,4.8mm>*{}**@{-},
     <0mm,0mm>*{};<-2mm,7mm>*{\overbrace{\ \ \ \ \ \ \ \ \ \ \ \ }}**@{},
     <0mm,0mm>*{};<-2mm,9mm>*{^{I_1}}**@{},
 <-0.6mm,-0.44mm>*{};<-8mm,-5mm>*{}**@{-},
 <-0.4mm,-0.7mm>*{};<-4.5mm,-5mm>*{}**@{-},
 <0mm,0mm>*{};<-1mm,-5mm>*{\ldots}**@{},
 <0.4mm,-0.7mm>*{};<4.5mm,-5mm>*{}**@{-},
 <0.6mm,-0.44mm>*{};<8mm,-5mm>*{}**@{-},
      <0mm,0mm>*{};<0mm,-7mm>*{\underbrace{\ \ \ \ \ \ \ \ \ \ \ \ \ \ \
      }}**@{},
      <0mm,0mm>*{};<0mm,-10.6mm>*{_{J_1}}**@{},
 <13mm,5mm>*{};<13mm,5mm>*{\circ}**@{},
 <12.6mm,5.44mm>*{};<5mm,10mm>*{}**@{-},
 <12.6mm,5.7mm>*{};<8.5mm,10mm>*{}**@{-},
 <13mm,5mm>*{};<13mm,10mm>*{\ldots}**@{},
 <13.4mm,5.7mm>*{};<16.5mm,10mm>*{}**@{-},
 <13.6mm,5.44mm>*{};<20mm,10mm>*{}**@{-},
      <13mm,5mm>*{};<13mm,12mm>*{\overbrace{\ \ \ \ \ \ \ \ \ \ \ \ \ \ }}**@{},
      <13mm,5mm>*{};<13mm,14mm>*{^{I_2}}**@{},
 <12.4mm,4.3mm>*{};<8mm,0mm>*{}**@{-},
 <12.6mm,4.3mm>*{};<12mm,0mm>*{\ldots}**@{},
 <13.4mm,4.5mm>*{};<16.5mm,0mm>*{}**@{-},
 <13.6mm,4.8mm>*{};<20mm,0mm>*{}**@{-},
     <13mm,5mm>*{};<14.3mm,-2mm>*{\underbrace{\ \ \ \ \ \ \ \ \ \ \ }}**@{},
     <13mm,5mm>*{};<14.3mm,-4.5mm>*{_{J_2}}**@{},
 \end{xy}}
$$
It was shown in \cite{Me1,Me2} that
representations $\HoLnB\rar \cE nd_V$ of the minimal resolution of
$\LnB$ in a graded vector space $V$ are in 1-1 correspondence with formal graded Poisson structures $\pi\in \cT_{poly}^{\geq 1}(V^*)$ on the dual vector space $V^*$ (viewed as a linear manifold) which vanish at the zero point in $V$, $\pi|_0=0$.

\subsection{\bf Quantizable odd Lie bialgebras} We define the properad $\LonB$ of {\em quantizable}\, odd Lie bialgebras as the quotient
of the  properad $\LnB$ by the ideal generated by the following element
\begin{equation}\label{equ:obstr_element}
\Ba{c}
\begin{xy}
 <0mm,-1.3mm>*{};<0mm,-3.5mm>*{}**@{-},
 <0.4mm,0.0mm>*{};<2.4mm,2.1mm>*{}**@{-},
 <-0.38mm,-0.2mm>*{};<-2.8mm,2.5mm>*{}**@{-},
<0mm,-0.8mm>*{\circ};<0mm,0.8mm>*{}**@{},
 <2.96mm,2.4mm>*{\circ};<2.45mm,2.35mm>*{}**@{},
 <2.4mm,2.8mm>*{};<0mm,5mm>*{}**@{-},
 <3.35mm,2.9mm>*{};<5.5mm,5mm>*{}**@{-},
<-2.8mm,2.5mm>*{};<0mm,5mm>*{\circ}**@{},
<-2.8mm,2.5mm>*{};<0mm,5mm>*{}**@{-},
<2.96mm,5mm>*{};<2.96mm,7.5mm>*{\circ}**@{},
<0.2mm,5.1mm>*{};<2.8mm,7.5mm>*{}**@{-},
<0.2mm,5.1mm>*{};<2.8mm,7.5mm>*{}**@{-},
<5.5mm,5mm>*{};<2.8mm,7.5mm>*{}**@{-},
<2.9mm,7.5mm>*{};<2.9mm,10.5mm>*{}**@{-},
    \end{xy}
\Ea\ \in \LnB.
\end{equation}
 The associated relation on Lie and coLie brackets looks like a higher genus odd analogue of the involutivity condition (\ref{1: involutivity condition})
 in the case of even Lie bialgebras. However, we prefer to use the adjective {\em quantizable}\, rather than {\em involutive} for odd Lie bialgebras satisfying (\ref{1: quantizability condition}) because that condition has a clear interpretation within the framework of the theory of deformation quantization, and its quantizability property becomes even more clear when one raises it to the level
 of  representations  of its minimal resolution $\HoLonB$.

\sip

An odd Lie bialgebra structure in a vector space $V$ can be understood as a pair
$$
\left(\xi\in \cT_{V^*}, \Phi\in \wedge^2 \cT_{V^*}\right)
$$
consisting of a degree 1 quadratic vector field $\xi$ (corresponding to the Lie cobracket $\Delta$ in $V$)
 and a linear Poisson structure $\Phi$ in $V^*$ (corresponding to the Lie bracket $[\ ,\ ]$ in $V$). All the (compatibility) equations for the algebraic operations $[\ ,\ ]$ and $\Delta$
get encoded into a single equation,
$$
\{\xi + \Phi, \xi + \Phi\}=0,
$$
where $\{\ ,\ \}$ stand for the standard Schouten bracket in the algebra $\cT_{poly}(V^*)$
of polyvector fields on $V^*$ (viewed as an affine manifold). Therefore, the sum $\xi + \Phi$ gives us a graded Poisson structure on $V^*$ and  one can talk about
its deformation quantization, that is, about
an associated Maurer-Cartan element $\Ga$ (deforming $\xi + \Phi$) in the Hochschild dg Lie algebra,
$$
C^\bu(\f_V,\f_V):= \bigoplus_{n\geq 0} \Hom(\f_V^{\ot n}, \f_V)
$$
where $\Hom(\f_V^{\ot n}, \f_V)$ stands for the vector space of polydifferential operators
on the graded commutative algebra $\f_V:=\odot^\bu V$ of polynomial functions on $V^*$. As the graded Poisson structure $\xi+\Phi$ is non-negatively graded, its deformation quantization must satisfy the condition
$$
\Ga\in \Hom(\K, \f_V) \oplus \Hom(\f_V,\f_V) \oplus \Hom(\f_V^{\ot 2},\f_V)
$$
with the corresponding splitting of $\Ga$ into a sum of three terms,
$$
\Ga=\Ga_0 + \Ga_1 + \Ga_2
$$
of degrees $2$, $1$, and $0$ respectively. The term $\Ga_2=\Ga_2(\Phi)$ has degree zero and hence can depend (universally) only on the Lie bracket $\Phi$. It makes $(\f_V, \star:=\Ga_2)$ into an associative non-commutative algebra, and up to  gauge equivalence, the algebra $(\f_V, \star)$ can always be identified
with the universal enveloping algebra of the Lie algebra $(V,[ \ ,\ ])$.
The operation $\Ga_1$ is a deformation of the differential on $\f_V$ induced by $\xi$. (Note that this latter undeformed differential squares to zero by the second Jacobi identity in the list (\ref{R for LieB}).)
The Maurer-Cartan equation for $\Ga$ states that $\Ga_1$ is a derivation with respect to the star product $\Ga_2$, and  squares to zero modulo the star product commutator with the (potential) obstruction $\Gamma_0$.
Now one can ask whether it is possible to find $\Ga$ as above such that the obstruction (or sometimes called "curvature") term $\Ga_0$ vanishes, and hence $\Ga_1^2=0$, so that $\Gamma_1$ is an honest differential.
Since the algebra $(\f_V, \star)$ is generated by $V$, any derivation with respect to the product $\star$ is uniquely determined by its values on $V$.
Let $\xi_\star$ be the unique derivation of the associative algebra $(\f_V, \star)$ that agrees with $\xi$ on $V$. This derivation is well defined since it annihilates the defining relations of the universal enveloping algebra by the third relation in (\ref{R for LieB}).
Then the derivation $\xi_\star^2 = \frac 1 2 [\xi_\star,\xi_\star]=0$ if and only if $\xi \simeq \begin{xy}
 <0mm,0.66mm>*{};<0mm,3mm>*{}**@{-},
 <0.39mm,-0.39mm>*{};<2.2mm,-2.2mm>*{}**@{-},
 <-0.35mm,-0.35mm>*{};<-2.2mm,-2.2mm>*{}**@{-},
 <0mm,0mm>*{\circ};<0mm,0mm>*{}**@{},
\end{xy}$ and $\Phi\simeq \begin{xy}
 <0mm,-0.55mm>*{};<0mm,-2.5mm>*{}**@{-},
 <0.5mm,0.5mm>*{};<2.2mm,2.2mm>*{}**@{-},
 <-0.48mm,0.48mm>*{};<-2.2mm,2.2mm>*{}**@{-},
 <0mm,0mm>*{\circ};<0mm,0mm>*{}**@{},
\end{xy}$ satisfy the extra compatibility condition (\ref{1: quantizability condition}) (see \cite{Me3}). Therefore, if $\xi$ and $\Phi$ come from a representation of $\LonB$ in $V$, then they admit a very simple deformation quantization in the form
$$
\Ga= \xi_\star + \Ga_2(\Phi),
$$
and this quantization makes sense even in the case when $V$ is {\em infinite}-dimensional.

\sip

If $\xi$ and $\Phi$ do not satisfy the extra compatibility condition (\ref{1: quantizability condition}), then their deformation quantization is possible only in {\em finite}\, dimensions, and
involve a non-zero ``curvature" term $\Ga_0$ which in turn involves graphs with closed paths of directed edges and is given explicitly in \cite{Me3} (in fact this argument proves non-existence of Kontsevich formality maps for {\em infinite}\, dimensional manifolds).

\sip

These considerations shall motivate our notation ``quantizable Lie bialgebras'' for odd Lie bialgebras satisfying the condition (\ref{1: quantizability condition}).

\sip

We shall construct below a minimal resolution $\HoLonB$ of the properad $\LonB$; its representations in a graded vector space $V$ give us so called {\em quantizable}\,
Poisson structures on $V$ which can be  deformation quantized via a trivial  (i.e.\ without using Drinfeld associators) perturbation even if $\dim V=\infty$ (see \cite{Wi2,B}); in finite dimensions there is a 1-1 correspondence between ordinary Poisson structure on $V$ and quantizable ones, but this correspondence is highly non-trivial --- it depends on the choice of a Drinfeld associator \cite{MW2}.

\sip

The properad Koszul dual to the properad $\LnB$ is the properad of odd Frobenius algebras (cf.\ \cite{V}).
A remarkable ``Koszul dual" meaning of the graph (\ref{1: quantizability condition}) was found by
Theo Johnson-Freyd in \cite{JF1} --- it controls the obstruction to the existence of a geometrically
meaningful homotopy odd Frobenius structure on the complex $\mathrm{Chains}_\bu(\R)$.

\bip

\bip

{\Large
\section{\bf A minimal resolution of $\LonB$}
}
\bip

\subsection{Oriented graph complexes and a Kontsevich-Shoikhet MC element}
Let $G_{n,l}^{or}$ be a set of connected graphs $\Ga$ with $n$ vertices and $l$ directed edges such that
(i) $\Ga$ has no {\em closed}\, directed paths of edges, and (ii)
some bijection from the set of edges $E(\Ga)$ to the set $[l]$ is fixed. There is
a natural right action of the group $\bS_l$ on the set $G^{or}_{n,l}$ by relabeling the
edges.

\sip

 Consider a graded vector space
$$
\mathsf{fGC}_2^{or}:= \prod_{n\geq 1, l\geq 0} \K \langle G_{n,l}^{or}\rangle \ot_{\bS_l}  \sgn_l [l +2(1-n)]
$$
It was shown in \cite{Wi2} that this vector space comes equipped with a Lie bracket $[\ ,\ ]$ (given, as
often in the theory of graph complexes, by substituting graphs into vertices of another graphs),
and that the degree $+1$ graph
$$\xy
 (0,0)*{\bullet}="a",
(6,0)*{\bu}="b",
\ar @{->} "a";"b" <0pt>
\endxy \in \mathsf{fGC}_2^{or}
$$
is a Maurer-Cartan element making $\mathsf{fGC}_2^{or}$ into a {\em differential}\, Lie algebra
with the differential given by
$$
\delta:=[\xy
 (0,0)*{\bullet}="a",
(6,0)*{\bu}="b",
\ar @{->} "a";"b" <0pt>
\endxy, \ ].
$$
 It was proven in \cite{Wi2} that the cohomology group $H^1(\fGC_2^{or})$ is one-dimensional and is spanned by the following graph
$$
\Upsilon_4:=\Ba{c}\xy
(0,0)*{\bullet}="1",
(-7,16)*{\bullet}="2",
(7,16)*{\bullet}="3",
(0,10)*{\bullet}="4",
\ar @{<-} "2";"4" <0pt>
\ar @{<-} "3";"4" <0pt>
\ar @{<-} "4";"1" <0pt>
\ar @{<-} "2";"1" <0pt>
\ar @{<-} "3";"1" <0pt>
\endxy\Ea
+
2
\Ba{c}\xy
(0,0)*{\bullet}="1",
(-6,6)*{\bullet}="2",
(6,10)*{\bullet}="3",
(0,16)*{\bullet}="4",
\ar @{<-} "4";"3" <0pt>
\ar @{<-} "4";"2" <0pt>
\ar @{<-} "3";"2" <0pt>
\ar @{<-} "2";"1" <0pt>
\ar @{<-} "3";"1" <0pt>
\endxy\Ea
+
 \Ba{c}\xy
(0,16)*{\bullet}="1",
(-7,0)*{\bullet}="2",
(7,0)*{\bullet}="3",
(0,6)*{\bullet}="4",
\ar @{->} "2";"4" <0pt>
\ar @{->} "3";"4" <0pt>
\ar @{->} "4";"1" <0pt>
\ar @{->} "2";"1" <0pt>
\ar @{->} "3";"1" <0pt>
\endxy\Ea.
$$
while the cohomology group $H^2(\mathsf{fGC}_2^{or},\delta_0)$ is also one-dimensional and is generated by a linear combination of graphs with four vertices (whose explicit form plays no role in this paper). This means that one can construct by induction a new Maurer-Cartan element (the integer subscript stand for the number of vertices)
$$
\Upsilon_{KS}= \xy
 (0,0)*{\bullet}="a",
(6,0)*{\bu}="b",
\ar @{->} "a";"b" <0pt>
\endxy  + \Upsilon_4
+ \Upsilon_6 + \Upsilon_8 + \ldots
$$
in the Lie algebra $\fGC_2^{or}$. Indeed, the Lie brackets in $\fGC_2^{or}$ has the property that
a commutator $[A,B]$ of a graph $A$ with $p$ vertices and a graph $B$ with $q$ vertices has $p+q-1$ vertices. Therefore, all the obstructions to extending the sum $ \xy
 (0,0)*{\bullet}="a",
(6,0)*{\bu}="b",
\ar @{->} "a";"b" <0pt>
\endxy  + \Upsilon_4 $ to a Maurer-Cartan element have $7$ or more vertices and hence do not hit the unique cohomology class
in $H^2(\mathsf{fGC}_2^{or},\delta)$. Up to gauge equivalence, this new MC element $\Upsilon_{KS}$ is the {\em only}\, non-trivial deformation of the standard MC element $\xy
 (0,0)*{\bullet}="a",
(6,0)*{\bu}="b",
\ar @{->} "a";"b" <0pt>
\endxy$. We call it the {\em Kontsevich-Shoikhet}\, element as it  was  introduced (via a different line of thought) by Boris Shoikhet in \cite{Sh} with a reference to an important contribution by Maxim Kontsevich via an informal communication.

\subsubsection{\bf A formal power series extension of $\fGC_2^{or}$}
Let $\hbar$
be a formal parameter of degree $0$ and let $\fGC_2^{or}[[\hbar]]$ be a topological vector space
of formal power series in $\hbar$ with coefficients in $\fGC_2^{or}$. This is naturally a topological Lie algebra in which the formal power series
$$
\Upsilon_{KS}^\hbar= \xy
 (0,0)*{\bullet}="a",
(6,0)*{\bu}="b",
\ar @{->} "a";"b" <0pt>
\endxy  + \hbar\Upsilon_4 +
+ \hbar^2 \Upsilon_6 + \hbar^3 \Upsilon_8 + \ldots
$$
is a Maurer-Cartan element.

\sip

\subsection{From the Kontsevich-Shoikhet element to a minimal resolution of $\LonB$}
Consider a (non-differential) free properad $\HoLonB$ generated by the following (skewsymmetric in outputs and symmetric in inputs) corollas  of degree $2-m$,
\Beq\label{2: generating corollas of LoB infty}
\resizebox{15mm}{!}{
\xy
(-9,-6)*{};
(0,0)*+{a}*\cir{}
**\dir{-};
(-5,-6)*{};
(0,0)*+{a}*\cir{}
**\dir{-};
(9,-6)*{};
(0,0)*+{a}*\cir{}
**\dir{-};
(5,-6)*{};
(0,0)*+{a}*\cir{}
**\dir{-};
(0,-6)*{\ldots};
(-10,-8)*{_1};
(-6,-8)*{_2};
(10,-8)*{_n};
(-9,6)*{};
(0,0)*+{a}*\cir{}
**\dir{-};
(-5,6)*{};
(0,0)*+{a}*\cir{}
**\dir{-};
(9,6)*{};
(0,0)*+{a}*\cir{}
**\dir{-};
(5,6)*{};
(0,0)*+{a}*\cir{}
**\dir{-};
(0,6)*{\ldots};
(-10,8)*{_1};
(-6,8)*{_2};
(10,8)*{_m};
\endxy}
=(-1)^{\sigma}
\resizebox{18mm}{!}{
\xy
(-9,-6)*{};
(0,0)*+{a}*\cir{}
**\dir{-};
(-5,-6)*{};
(0,0)*+{a}*\cir{}
**\dir{-};
(9,-6)*{};
(0,0)*+{a}*\cir{}
**\dir{-};
(5,-6)*{};
(0,0)*+{a}*\cir{}
**\dir{-};
(0,-6)*{\ldots};
(-12,-8)*{_{\tau(1)}};
(-6,-8)*{_{\tau(2)}};
(12,-8)*{_{\tau(n)}};
(-9,6)*{};
(0,0)*+{a}*\cir{}
**\dir{-};
(-5,6)*{};
(0,0)*+{a}*\cir{}
**\dir{-};
(9,6)*{};
(0,0)*+{a}*\cir{}
**\dir{-};
(5,6)*{};
(0,0)*+{a}*\cir{}
**\dir{-};
(0,6)*{\ldots};
(-12,8)*{_{\sigma(1)}};
(-6,8)*{_{\sigma(2)}};
(12,8)*{_{\sigma(m)}};
\endxy}\ \ \ \forall \sigma\in \bS_m, \forall \tau\in \bS_n,
\Eeq
where $m+n+ a\geq 3$, $m\geq 1$, $n\geq 1$, $a\geq 0$. Let $\wHoLonB$ be the genus completion of
$\HoLonB$.

\subsubsection{\bf Lemma} {\em The Lie algebra $\fGC_2^{or}[[\hbar]]$ acts (from the right) on the properad $\wHoLonB$ by continuous derivations, that is, there is a morphism of Lie algebras
$$
\Ba{rccc}
F: & \fGC_2^{or}[[\hbar]] & \lon & \Der(\wHoLonB)\\
   & \hbar^k \Ga &\lon & F(\hbar^k \Ga)
\Ea
$$
where the derivation $F(\hbar^k \Ga)$ is given on the generators as follows
$$
F(\hbar^k \Ga) \cdot \resizebox{15mm}{!}{
\xy
(-9,-6)*{};
(0,0)*+{a}*\cir{}
**\dir{-};
(-5,-6)*{};
(0,0)*+{a}*\cir{}
**\dir{-};
(9,-6)*{};
(0,0)*+{a}*\cir{}
**\dir{-};
(5,-6)*{};
(0,0)*+{a}*\cir{}
**\dir{-};
(0,-6)*{\ldots};
(-10,-8)*{_1};
(-6,-8)*{_2};
(10,-8)*{_n};
(-9,6)*{};
(0,0)*+{a}*\cir{}
**\dir{-};
(-5,6)*{};
(0,0)*+{a}*\cir{}
**\dir{-};
(9,6)*{};
(0,0)*+{a}*\cir{}
**\dir{-};
(5,6)*{};
(0,0)*+{a}*\cir{}
**\dir{-};
(0,6)*{\ldots};
(-10,8)*{_1};
(-6,8)*{_2};
(10,8)*{_m};
\endxy}
:=
\left\{\Ba{ll} \displaystyle
\sum_{m,n}
 \sum_{a=a_1+\ldots+ a_{\#V(\Ga)} +k }
    \overbrace{
 \underbrace{\Ba{c}\resizebox{10mm}{!}  { \xy
(0,4.5)*+{...},
(0,-4.5)*+{...},
(0,0)*+{\Ga}="o",
(-5,6)*{}="1",
(-3,6)*{}="2",
(3,6)*{}="3",
(5,6)*{}="4",
(-3,-6)*{}="5",
(3,-6)*{}="6",
(5,-6)*{}="7",
(-5,-6)*{}="8",
\ar @{-} "o";"1" <0pt>
\ar @{-} "o";"2" <0pt>
\ar @{-} "o";"3" <0pt>
\ar @{-} "o";"4" <0pt>
\ar @{-} "o";"5" <0pt>
\ar @{-} "o";"6" <0pt>
\ar @{-} "o";"7" <0pt>
\ar @{-} "o";"8" <0pt>
\endxy}\Ea
 }_{n\times}
 }^{m\times} & \ \ \ \ \ \ \forall\ \Ga\in \fGC_2^{or}, \ \ \ \ \forall\ k\in [0,1,\ldots, a]\\
 0 & \ \ \ \ \ \  \forall\ \Ga\in \fGC_2^{or},\ \ \ \ \forall\ k>a,
 \Ea
 \right.
$$
where the first sum is taking over to attach $m$ output legs and $n$ input legs to the vertices
of the graph $\Ga$, and the second sum is taken over all possible ways to decorate the vertices of $\Ga$ with non-negative integers $a_1,\ldots,a_{\# V(\Ga)}$ such they sum to $a-k$.
}

\mip

{ Proof} is identical to the proofs of similar statements (Theorems 1.2.1 and 1.2.2) in \cite{MW}.

\subsubsection{\bf Corollary} {\em The completed free properad $\wHoLonB$ comes equipped with a differential
$\delta_\diamond:=F(\Upsilon_{KS}^\hbar)$. The differential $\delta$ restricts to a differential in the  free
properad $\HoLonB$}.

\begin{proof} When applied to any generator of $\wHoLonB$ the differential $\delta$ gives always a {\em finite}\, sum of graphs. It follows that it is well defined in $\HoLonB$ as well.
\end{proof}

\bip

There is an injection of dg free properads
$$
(\HoLnB, \delta) \lon (\HoLonB, \delta_\diamond)
$$
given on generators by
$$
\resizebox{14mm}{!}{\begin{xy}
 <0mm,0mm>*{\circ};<0mm,0mm>*{}**@{},
 <-0.6mm,0.44mm>*{};<-8mm,5mm>*{}**@{-},
 <-0.4mm,0.7mm>*{};<-4.5mm,5mm>*{}**@{-},
 <0mm,0mm>*{};<-1mm,5mm>*{\ldots}**@{},
 <0.4mm,0.7mm>*{};<4.5mm,5mm>*{}**@{-},
 <0.6mm,0.44mm>*{};<8mm,5mm>*{}**@{-},
   <0mm,0mm>*{};<-8.5mm,5.5mm>*{^1}**@{},
   <0mm,0mm>*{};<-5mm,5.5mm>*{^2}**@{},
   <0mm,0mm>*{};<4.5mm,5.5mm>*{^{m\hspace{-0.5mm}-\hspace{-0.5mm}1}}**@{},
   <0mm,0mm>*{};<9.0mm,5.5mm>*{^m}**@{},
 <-0.6mm,-0.44mm>*{};<-8mm,-5mm>*{}**@{-},
 <-0.4mm,-0.7mm>*{};<-4.5mm,-5mm>*{}**@{-},
 <0mm,0mm>*{};<-1mm,-5mm>*{\ldots}**@{},
 <0.4mm,-0.7mm>*{};<4.5mm,-5mm>*{}**@{-},
 <0.6mm,-0.44mm>*{};<8mm,-5mm>*{}**@{-},
   <0mm,0mm>*{};<-8.5mm,-6.9mm>*{^1}**@{},
   <0mm,0mm>*{};<-5mm,-6.9mm>*{^2}**@{},
   <0mm,0mm>*{};<4.5mm,-6.9mm>*{^{n\hspace{-0.5mm}-\hspace{-0.5mm}1}}**@{},
   <0mm,0mm>*{};<9.0mm,-6.9mm>*{^n}**@{},
 \end{xy}}
  \lon
 \resizebox{15mm}{!}{
\xy
(-9,-6)*{};
(0,0)*+{0}*\cir{}
**\dir{-};
(-5,-6)*{};
(0,0)*+{0}*\cir{}
**\dir{-};
(9,-6)*{};
(0,0)*+{0}*\cir{}
**\dir{-};
(5,-6)*{};
(0,0)*+{0}*\cir{}
**\dir{-};
(0,-6)*{\ldots};
(-10,-8)*{_1};
(-6,-8)*{_2};
(10,-8)*{_n};
(-9,6)*{};
(0,0)*+{0}*\cir{}
**\dir{-};
(-5,6)*{};
(0,0)*+{0}*\cir{}
**\dir{-};
(9,6)*{};
(0,0)*+{0}*\cir{}
**\dir{-};
(5,6)*{};
(0,0)*+{0}*\cir{}
**\dir{-};
(0,6)*{\ldots};
(-10,8)*{_1};
(-6,8)*{_2};
(10,8)*{_m};
\endxy}
$$
Identifying from now on weight zero generators of $\HoLonB$ with  generators of $\HoLnB$, we may
write
$$
\delta_\diamond
\xy
(0,5)*{};
(0,0)*+{_1}*\cir{}
**\dir{-};
(0,-5)*{};
(0,0)*+{_1}*\cir{}
**\dir{-};
\endxy
=\Ba{c}
\begin{xy}
 <0mm,-1.3mm>*{};<0mm,-3.5mm>*{}**@{-},
 <0.4mm,0.0mm>*{};<2.4mm,2.1mm>*{}**@{-},
 <-0.38mm,-0.2mm>*{};<-2.8mm,2.5mm>*{}**@{-},
<0mm,-0.8mm>*{\circ};<0mm,0.8mm>*{}**@{},
 <2.96mm,2.4mm>*{\circ};<2.45mm,2.35mm>*{}**@{},
 <2.4mm,2.8mm>*{};<0mm,5mm>*{}**@{-},
 <3.35mm,2.9mm>*{};<5.5mm,5mm>*{}**@{-},
<-2.8mm,2.5mm>*{};<0mm,5mm>*{\circ}**@{},
<-2.8mm,2.5mm>*{};<0mm,5mm>*{}**@{-},
<2.96mm,5mm>*{};<2.96mm,7.5mm>*{\circ}**@{},
<0.2mm,5.1mm>*{};<2.8mm,7.5mm>*{}**@{-},
<0.2mm,5.1mm>*{};<2.8mm,7.5mm>*{}**@{-},
<5.5mm,5mm>*{};<2.8mm,7.5mm>*{}**@{-},
<2.9mm,7.5mm>*{};<2.9mm,10.5mm>*{}**@{-},
    \end{xy}
    \vspace{4mm}
\Ea
$$
and hence conclude that there is a natural morphism of dg properads
$$
\pi: (\HoLonB,\delta_\diamond) \lon (\LonB,0)
$$
Our main result in this paper is the following theorem.

\subsection{\bf Main Theorem} {\em The map $\pi$ is a quasi-isomorphism, i.e.\ $\HoLonB$ is a minimal resolution of $\LonB$.}

\begin{proof}
Let $\cP$ be a dg properad
generated by a degree $1$ corollas
$\begin{xy}
 <0mm,-0.55mm>*{};<0mm,-3mm>*{}**@{-},
 <0mm,0.5mm>*{};<0mm,3mm>*{}**@{-},
 <0mm,0mm>*{\bu};<0mm,0mm>*{}**@{},
 \end{xy}$\, and
 $\Ba{c}\begin{xy}
 <0mm,0.66mm>*{};<0mm,3mm>*{}**@{-},
 <0.39mm,-0.39mm>*{};<2.2mm,-2.2mm>*{}**@{-},
 <-0.35mm,-0.35mm>*{};<-2.2mm,-2.2mm>*{}**@{-},
 <0mm,0mm>*{\circ};<0mm,0mm>*{}**@{},
   <0.39mm,-0.39mm>*{};<2.9mm,-4mm>*{^2}**@{},
   <-0.35mm,-0.35mm>*{};<-2.8mm,-4mm>*{^1}**@{},
\end{xy}=
\begin{xy}
 <0mm,0.66mm>*{};<0mm,3mm>*{}**@{-},
 <0.39mm,-0.39mm>*{};<2.2mm,-2.2mm>*{}**@{-},
 <-0.35mm,-0.35mm>*{};<-2.2mm,-2.2mm>*{}**@{-},
 <0mm,0mm>*{\circ};<0mm,0mm>*{}**@{},
   <0.39mm,-0.39mm>*{};<2.9mm,-4mm>*{^1}**@{},
   <-0.35mm,-0.35mm>*{};<-2.8mm,-4mm>*{^2}**@{},
\end{xy}\Ea$,
 and a degree zero corolla,
$
\Ba{c}
 \begin{xy}
 <0mm,-0.55mm>*{};<0mm,-2.5mm>*{}**@{-},
 <0.5mm,0.5mm>*{};<2.2mm,2.2mm>*{}**@{-},
 <-0.48mm,0.48mm>*{};<-2.2mm,2.2mm>*{}**@{-},
 <0mm,0mm>*{\circ};<0mm,0mm>*{}**@{},
 <0.5mm,0.5mm>*{};<2.7mm,2.8mm>*{^2}**@{},
 <-0.48mm,0.48mm>*{};<-2.7mm,2.8mm>*{^1}**@{},
 \end{xy}
=-
\begin{xy}
 <0mm,-0.55mm>*{};<0mm,-2.5mm>*{}**@{-},
 <0.5mm,0.5mm>*{};<2.2mm,2.2mm>*{}**@{-},
 <-0.48mm,0.48mm>*{};<-2.2mm,2.2mm>*{}**@{-},
 <0mm,0mm>*{\circ};<0mm,0mm>*{}**@{},
 <0.5mm,0.5mm>*{};<2.7mm,2.8mm>*{^1}**@{},
 <-0.48mm,0.48mm>*{};<-2.7mm,2.8mm>*{^2}**@{},
 \end{xy}\Ea
 $ and
 modulo relations
 \Beq\label{2: relation in P}
\xy
(0,-1.9)*{\bu}="0",
 (0,1.9)*{\bu}="1",
(0,-5)*{}="d",
(0,5)*{}="u",
\ar @{-} "0";"u" <0pt>
\ar @{-} "0";"1" <0pt>
\ar @{-} "1";"d" <0pt>
\endxy=0, \ \ \ \
\xy
(0,-1.9)*{\circ}="0",
 (0,1.9)*{\bu}="1",
(-2.5,-5)*{}="d1",
(2.5,-5)*{}="d2",
(0,5)*{}="u",
\ar @{-} "1";"u" <0pt>
\ar @{-} "0";"1" <0pt>
\ar @{-} "0";"d1" <0pt>
\ar @{-} "0";"d2" <0pt>
\endxy  +
\xy
(-2.5,-1.9)*{\bu}="0",
 (0,1.9)*{\circ}="1",
(-2.5,-1.9)*{}="d1",
(2.5,-1.9)*{}="d2",
(-2.5,-5)*{}="d",
(0,5)*{}="u",
\ar @{-} "1";"u" <0pt>
\ar @{-} "0";"1" <0pt>
\ar @{-} "0";"d" <0pt>
\ar @{-} "1";"d1" <0pt>
\ar @{-} "1";"d2" <0pt>
\endxy
  +
\xy
(2.5,-1.9)*{\bu}="0",
 (0,1.9)*{\circ}="1",
(-2.5,-1.9)*{}="d1",
(2.5,-1.9)*{}="d2",
(2.5,-5)*{}="d",
(0,5)*{}="u",
\ar @{-} "1";"u" <0pt>
\ar @{-} "0";"1" <0pt>
\ar @{-} "0";"d" <0pt>
\ar @{-} "1";"d1" <0pt>
\ar @{-} "1";"d2" <0pt>
\endxy
=0\ , \ \ \ \
\xy
(0,1.9)*{\circ}="0",
 (0,-1.9)*{\bu}="1",
(-2.5,5)*{}="d1",
(2.5,5)*{}="d2",
(0,-5)*{}="u",
\ar @{-} "1";"u" <0pt>
\ar @{-} "0";"1" <0pt>
\ar @{-} "0";"d1" <0pt>
\ar @{-} "0";"d2" <0pt>
\endxy  -
\xy
(-2.5,1.9)*{\bu}="0",
 (0,-1.9)*{\circ}="1",
(-2.5,1.9)*{}="d1",
(2.5,1.9)*{}="d2",
(-2.5,5)*{}="d",
(0,-5)*{}="u",
\ar @{-} "1";"u" <0pt>
\ar @{-} "0";"1" <0pt>
\ar @{-} "0";"d" <0pt>
\ar @{-} "1";"d1" <0pt>
\ar @{-} "1";"d2" <0pt>
\endxy
  -
\xy
(2.5,1.9)*{\bu}="0",
 (0,-1.9)*{\circ}="1",
(-2.5,1.9)*{}="d1",
(2.5,1.9)*{}="d2",
(2.5,5)*{}="d",
(0,-5)*{}="u",
\ar @{-} "1";"u" <0pt>
\ar @{-} "0";"1" <0pt>
\ar @{-} "0";"d" <0pt>
\ar @{-} "1";"d1" <0pt>
\ar @{-} "1";"d2" <0pt>
\endxy
=0.
\Eeq and the  three relations in (\ref{R for
 LieB}).
 The differential in $\cP$ is given on the generators by
\Beq\label{2: differential in P}
d\, \begin{xy}
 <0mm,0.66mm>*{};<0mm,3mm>*{}**@{-},
 <0.39mm,-0.39mm>*{};<2.2mm,-2.2mm>*{}**@{-},
 <-0.35mm,-0.35mm>*{};<-2.2mm,-2.2mm>*{}**@{-},
 <0mm,0mm>*{\circ};<0mm,0mm>*{}**@{},
\end{xy}=0, \ \ \ \
d\, \begin{xy}
 <0mm,-0.55mm>*{};<0mm,-2.5mm>*{}**@{-},
 <0.5mm,0.5mm>*{};<2.2mm,2.2mm>*{}**@{-},
 <-0.48mm,0.48mm>*{};<-2.2mm,2.2mm>*{}**@{-},
 <0mm,0mm>*{\circ};<0mm,0mm>*{}**@{},
 \end{xy} =0, \ \ \ \
d\, \begin{xy}
 <0mm,-0.55mm>*{};<0mm,-3mm>*{}**@{-},
 <0mm,0.5mm>*{};<0mm,3mm>*{}**@{-},
 <0mm,0mm>*{\bu};<0mm,0mm>*{}**@{},
 \end{xy}=
 \Ba{c}
\begin{xy}
 <0mm,-1.3mm>*{};<0mm,-3.5mm>*{}**@{-},
 <0.4mm,0.0mm>*{};<2.4mm,2.1mm>*{}**@{-},
 <-0.38mm,-0.2mm>*{};<-2.8mm,2.5mm>*{}**@{-},
<0mm,-0.8mm>*{\circ};<0mm,0.8mm>*{}**@{},
 <2.96mm,2.4mm>*{\circ};<2.45mm,2.35mm>*{}**@{},
 <2.4mm,2.8mm>*{};<0mm,5mm>*{}**@{-},
 <3.35mm,2.9mm>*{};<5.5mm,5mm>*{}**@{-},
<-2.8mm,2.5mm>*{};<0mm,5mm>*{\circ}**@{},
<-2.8mm,2.5mm>*{};<0mm,5mm>*{}**@{-},
<2.96mm,5mm>*{};<2.96mm,7.5mm>*{\circ}**@{},
<0.2mm,5.1mm>*{};<2.8mm,7.5mm>*{}**@{-},
<0.2mm,5.1mm>*{};<2.8mm,7.5mm>*{}**@{-},
<5.5mm,5mm>*{};<2.8mm,7.5mm>*{}**@{-},
<2.9mm,7.5mm>*{};<2.9mm,10.5mm>*{}**@{-},
    \end{xy}
    \vspace{4mm}
\Ea
\Eeq
{\sc Claim I}: {\em The surjective morphism of
dg properads,
\Beq\label{2: nu quasi-iso to P}
\nu: \HoLonB \lon \cP,
\Eeq
which sends all generators to zero except for the following ones
\Beq\label{2: map from P to P}
\nu\left( \xy
(-3,-4)*{};
(0,0)*+{_0}*\cir{}
**\dir{-};
(3,-4)*{};
(0,0)*+{_0}*\cir{}
**\dir{-};
(0,5)*{};
(0,0)*+{_0}*\cir{}
**\dir{-};
\endxy\right)
=\begin{xy}
 <0mm,0.66mm>*{};<0mm,4mm>*{}**@{-},
 <0.39mm,-0.39mm>*{};<2.2mm,-3.2mm>*{}**@{-},
 <-0.35mm,-0.35mm>*{};<-2.2mm,-3.2mm>*{}**@{-},
 <0mm,0mm>*{\circ};<0mm,0mm>*{}**@{},
\end{xy}\ ,\ \ \ \ \ \
\nu\left( \xy
(-3,4)*{};
(0,0)*+{_0}*\cir{}
**\dir{-};
(3,4)*{};
(0,0)*+{_0}*\cir{}
**\dir{-};
(0,-5)*{};
(0,0)*+{_0}*\cir{}
**\dir{-};
\endxy\right)
=\begin{xy}
 <0mm,-0.66mm>*{};<0mm,-4mm>*{}**@{-},
 <0.4mm,0.4mm>*{};<2.2mm,3.2mm>*{}**@{-},
 <-0.4mm,0.4mm>*{};<-2.2mm,3.2mm>*{}**@{-},
 <0mm,-0.1mm>*{\circ};<0mm,0mm>*{}**@{},
\end{xy}\ \ , \ \
\nu\left(
\xy
(0,5)*{};
(0,0)*+{_1}*\cir{}
**\dir{-};
(0,-5)*{};
(0,0)*+{_1}*\cir{}
**\dir{-};
\endxy\right)=\begin{xy}
 <0mm,-0.55mm>*{};<0mm,-3mm>*{}**@{-},
 <0mm,0.5mm>*{};<0mm,3mm>*{}**@{-},
 <0mm,0mm>*{\bu};<0mm,0mm>*{}**@{},
 \end{xy}
\Eeq
is a quasi-isomorphism.}

\sip

The proof of this claim is identical to the proof of Theorem 2.7.1 in \cite{CMW} so that we can omit
the details.

\bip

The proof of the Main Theorem will be completed once we show the following

\sip
{\sc Claim II:} {\em The natural map
\Beq\label{3: Claim 2 map mu}
\mu: (\cP, d) \lon (\LonB,0)
\Eeq
is a quasi-isomorphism.}

\bip

Let us define a new homological grading in the properad $\cP$ by assigning to the generator
$\begin{xy}
 <0mm,-0.55mm>*{};<0mm,-3mm>*{}**@{-},
 <0mm,0.5mm>*{};<0mm,3mm>*{}**@{-},
 <0mm,0mm>*{\bu};<0mm,0mm>*{}**@{},
 \end{xy}$ degree $-1$ and to the remaining generators the degree zero; to avoid confusion with the original grading let us call this new grading $s$-{\em grading}. Then  Claim II is proven once
 we show that  the cohomology $H(\cP)$ of $\cP$ is concentrated in $s$-degree zero.

  \sip

  Consider a path filtration  \cite{Ko,MaVo} of the dg properad $\cP$. The associated graded $\gr\cP$ can be identified with dg properad generated by the same corollas $\begin{xy}
 <0mm,-0.55mm>*{};<0mm,-3mm>*{}**@{-},
 <0mm,0.5mm>*{};<0mm,3mm>*{}**@{-},
 <0mm,0mm>*{\bu};<0mm,0mm>*{}**@{},
 \end{xy}$,
 $\Ba{c}\begin{xy}
 <0mm,0.66mm>*{};<0mm,3mm>*{}**@{-},
 <0.39mm,-0.39mm>*{};<2.2mm,-2.2mm>*{}**@{-},
 <-0.35mm,-0.35mm>*{};<-2.2mm,-2.2mm>*{}**@{-},
 <0mm,0mm>*{\circ};<0mm,0mm>*{}**@{},
   <0.39mm,-0.39mm>*{};<2.9mm,-4mm>*{^2}**@{},
   <-0.35mm,-0.35mm>*{};<-2.8mm,-4mm>*{^1}**@{},
\end{xy}=
\begin{xy}
 <0mm,0.66mm>*{};<0mm,3mm>*{}**@{-},
 <0.39mm,-0.39mm>*{};<2.2mm,-2.2mm>*{}**@{-},
 <-0.35mm,-0.35mm>*{};<-2.2mm,-2.2mm>*{}**@{-},
 <0mm,0mm>*{\circ};<0mm,0mm>*{}**@{},
   <0.39mm,-0.39mm>*{};<2.9mm,-4mm>*{^1}**@{},
   <-0.35mm,-0.35mm>*{};<-2.8mm,-4mm>*{^2}**@{},
\end{xy}\Ea$
 and
$
\Ba{c}
 \begin{xy}
 <0mm,-0.55mm>*{};<0mm,-2.5mm>*{}**@{-},
 <0.5mm,0.5mm>*{};<2.2mm,2.2mm>*{}**@{-},
 <-0.48mm,0.48mm>*{};<-2.2mm,2.2mm>*{}**@{-},
 <0mm,0mm>*{\circ};<0mm,0mm>*{}**@{},
 <0.5mm,0.5mm>*{};<2.7mm,2.8mm>*{^2}**@{},
 <-0.48mm,0.48mm>*{};<-2.7mm,2.8mm>*{^1}**@{},
 \end{xy}
=-
\begin{xy}
 <0mm,-0.55mm>*{};<0mm,-2.5mm>*{}**@{-},
 <0.5mm,0.5mm>*{};<2.2mm,2.2mm>*{}**@{-},
 <-0.48mm,0.48mm>*{};<-2.2mm,2.2mm>*{}**@{-},
 <0mm,0mm>*{\circ};<0mm,0mm>*{}**@{},
 <0.5mm,0.5mm>*{};<2.7mm,2.8mm>*{^1}**@{},
 <-0.48mm,0.48mm>*{};<-2.7mm,2.8mm>*{^2}**@{},
 \end{xy}\Ea
 $
which are subject to the relations (\ref{2: relation in P}), the first two relation in
(\ref{R for LieB}) and the following one
$$
\begin{xy}
 <0mm,2.47mm>*{};<0mm,0.12mm>*{}**@{-},
 <0.5mm,3.5mm>*{};<2.2mm,5.2mm>*{}**@{-},
 <-0.48mm,3.48mm>*{};<-2.2mm,5.2mm>*{}**@{-},
 <0mm,3mm>*{\circ};<0mm,3mm>*{}**@{},
  <0mm,-0.8mm>*{\circ};<0mm,-0.8mm>*{}**@{},
<-0.39mm,-1.2mm>*{};<-2.2mm,-3.5mm>*{}**@{-},
 <0.39mm,-1.2mm>*{};<2.2mm,-3.5mm>*{}**@{-},
     <0.5mm,3.5mm>*{};<2.8mm,5.7mm>*{^2}**@{},
     <-0.48mm,3.48mm>*{};<-2.8mm,5.7mm>*{^1}**@{},
   <0mm,-0.8mm>*{};<-2.7mm,-5.2mm>*{^1}**@{},
   <0mm,-0.8mm>*{};<2.7mm,-5.2mm>*{^2}**@{},
\end{xy}=0
$$
The differential in $\gr\cP$ is given by the original formula (\ref{2: differential in P}).
The {\sc Claim II} is proven once it is shown that the cohomology of the properad $\gr\cP$ is concentrated
in $s$-degree zero, or equivalently, the cohomology of dg prop $\cU(\gr\cP)$ generated
by this properad is concentrated in $s$-degree zero (as the universal enveloping functor
$\cU$ from the category of properads to the category of props is exact).

\sip

Consider a free prop $\cF ree\langle E \rangle$ generated by an $\bS$-bimodule
$E=\{E(m,n)\}$ with all $E(m,n)=0$ except the following ones,

\Beqrn
E(1,1) &=& \K[-1]=\mathrm{span}\left\langle \begin{xy}
 <0mm,-0.55mm>*{};<0mm,-3mm>*{}**@{-},
 <0mm,0.5mm>*{};<0mm,3mm>*{}**@{-},
 <0mm,0mm>*{\bu};<0mm,0mm>*{}**@{},
 \end{xy}   \right\rangle \\
 E(2,1) &=& sgn_2 =\mathrm{span}\left\langle
 \Ba{c}
 \begin{xy}
 <0mm,-0.55mm>*{};<0mm,-2.5mm>*{}**@{-},
 <0.5mm,0.5mm>*{};<2.2mm,2.2mm>*{}**@{-},
 <-0.48mm,0.48mm>*{};<-2.2mm,2.2mm>*{}**@{-},
 <0mm,0mm>*{\circ};<0mm,0mm>*{}**@{},
 <0.5mm,0.5mm>*{};<2.7mm,2.8mm>*{^2}**@{},
 <-0.48mm,0.48mm>*{};<-2.7mm,2.8mm>*{^1}**@{},
 \end{xy}
=-
\begin{xy}
 <0mm,-0.55mm>*{};<0mm,-2.5mm>*{}**@{-},
 <0.5mm,0.5mm>*{};<2.2mm,2.2mm>*{}**@{-},
 <-0.48mm,0.48mm>*{};<-2.2mm,2.2mm>*{}**@{-},
 <0mm,0mm>*{\circ};<0mm,0mm>*{}**@{},
 <0.5mm,0.5mm>*{};<2.7mm,2.8mm>*{^1}**@{},
 <-0.48mm,0.48mm>*{};<-2.7mm,2.8mm>*{^2}**@{},
 \end{xy}\Ea\right\rangle\\
 E(1,2) &=& \K[\bS_2][-1] =\mathrm{span}\left\langle
 \Ba{c}\begin{xy}
 <0mm,0.66mm>*{};<0mm,3mm>*{}**@{-},
 <0.39mm,-0.39mm>*{};<2.2mm,-2.2mm>*{}**@{-},
 <-0.35mm,-0.35mm>*{};<-2.2mm,-2.2mm>*{}**@{-},
 <0mm,0mm>*{\bu};<0mm,0mm>*{}**@{},
   <0.39mm,-0.39mm>*{};<2.9mm,-4mm>*{^2}**@{},
   <-0.35mm,-0.35mm>*{};<-2.8mm,-4mm>*{^1}**@{},
\end{xy}\neq
\begin{xy}
 <0mm,0.66mm>*{};<0mm,3mm>*{}**@{-},
 <0.39mm,-0.39mm>*{};<2.2mm,-2.2mm>*{}**@{-},
 <-0.35mm,-0.35mm>*{};<-2.2mm,-2.2mm>*{}**@{-},
 <0mm,0mm>*{\bu};<0mm,0mm>*{}**@{},
   <0.39mm,-0.39mm>*{};<2.9mm,-4mm>*{^1}**@{},
   <-0.35mm,-0.35mm>*{};<-2.8mm,-4mm>*{^2}**@{},
\end{xy}\Ea
\right\rangle
\Eeqrn
We assign to the above generators $s$-degrees $-1$, $0$ and $0$ respectively.

\sip

Define next a dg prop $\cA$ as the quotient of the above free prop $\cF ree\langle E\rangle$
by the ideal generated by the relations
$$
\xy
(0,-1.9)*{\bu}="0",
 (0,1.9)*{\bu}="1",
(0,-5)*{}="d",
(0,5)*{}="u",
\ar @{-} "0";"u" <0pt>
\ar @{-} "0";"1" <0pt>
\ar @{-} "1";"d" <0pt>
\endxy=0, \ \ \ \
\xy
(0,-1.9)*{\bu}="0",
 (0,1.9)*{\bu}="1",
(-2.5,-5)*{}="d1",
(2.5,-5)*{}="d2",
(0,5)*{}="u",
\ar @{-} "1";"u" <0pt>
\ar @{-} "0";"1" <0pt>
\ar @{-} "0";"d1" <0pt>
\ar @{-} "0";"d2" <0pt>
\endxy  +
\xy
(-2.5,-1.9)*{\bu}="0",
 (0,1.9)*{\bu}="1",
(-2.5,-1.9)*{}="d1",
(2.5,-1.9)*{}="d2",
(-2.5,-5)*{}="d",
(0,5)*{}="u",
\ar @{-} "1";"u" <0pt>
\ar @{-} "0";"1" <0pt>
\ar @{-} "0";"d" <0pt>
\ar @{-} "1";"d1" <0pt>
\ar @{-} "1";"d2" <0pt>
\endxy
  +
\xy
(2.5,-1.9)*{\bu}="0",
 (0,1.9)*{\bu}="1",
(-2.5,-1.9)*{}="d1",
(2.5,-1.9)*{}="d2",
(2.5,-5)*{}="d",
(0,5)*{}="u",
\ar @{-} "1";"u" <0pt>
\ar @{-} "0";"1" <0pt>
\ar @{-} "0";"d" <0pt>
\ar @{-} "1";"d1" <0pt>
\ar @{-} "1";"d2" <0pt>
\endxy
=0\ , \ \ \ \
\xy
(0,1.9)*{\circ}="0",
 (0,-1.9)*{\bu}="1",
(-2.5,5)*{}="d1",
(2.5,5)*{}="d2",
(0,-5)*{}="u",
\ar @{-} "1";"u" <0pt>
\ar @{-} "0";"1" <0pt>
\ar @{-} "0";"d1" <0pt>
\ar @{-} "0";"d2" <0pt>
\endxy  -
\xy
(-2.5,1.9)*{\bu}="0",
 (0,-1.9)*{\circ}="1",
(-2.5,1.9)*{}="d1",
(2.5,1.9)*{}="d2",
(-2.5,5)*{}="d",
(0,-5)*{}="u",
\ar @{-} "1";"u" <0pt>
\ar @{-} "0";"1" <0pt>
\ar @{-} "0";"d" <0pt>
\ar @{-} "1";"d1" <0pt>
\ar @{-} "1";"d2" <0pt>
\endxy
  -
\xy
(2.5,1.9)*{\bu}="0",
 (0,-1.9)*{\circ}="1",
(-2.5,1.9)*{}="d1",
(2.5,1.9)*{}="d2",
(2.5,5)*{}="d",
(0,-5)*{}="u",
\ar @{-} "1";"u" <0pt>
\ar @{-} "0";"1" <0pt>
\ar @{-} "0";"d" <0pt>
\ar @{-} "1";"d1" <0pt>
\ar @{-} "1";"d2" <0pt>
\endxy
=0
$$
and
$$
\begin{xy}
 <0mm,0mm>*{\circ};<0mm,0mm>*{}**@{},
 <0mm,-0.49mm>*{};<0mm,-3.0mm>*{}**@{-},
 <0.49mm,0.49mm>*{};<1.9mm,1.9mm>*{}**@{-},
 <-0.5mm,0.5mm>*{};<-1.9mm,1.9mm>*{}**@{-},
 <-2.3mm,2.3mm>*{\circ};<-2.3mm,2.3mm>*{}**@{},
 <-1.8mm,2.8mm>*{};<0mm,4.9mm>*{}**@{-},
 <-2.8mm,2.9mm>*{};<-4.6mm,4.9mm>*{}**@{-},
   <0.49mm,0.49mm>*{};<2.7mm,2.3mm>*{^3}**@{},
   <-1.8mm,2.8mm>*{};<0.4mm,5.3mm>*{^2}**@{},
   <-2.8mm,2.9mm>*{};<-5.1mm,5.3mm>*{^1}**@{},
 \end{xy}
\ + \
\begin{xy}
 <0mm,0mm>*{\circ};<0mm,0mm>*{}**@{},
 <0mm,-0.49mm>*{};<0mm,-3.0mm>*{}**@{-},
 <0.49mm,0.49mm>*{};<1.9mm,1.9mm>*{}**@{-},
 <-0.5mm,0.5mm>*{};<-1.9mm,1.9mm>*{}**@{-},
 <-2.3mm,2.3mm>*{\circ};<-2.3mm,2.3mm>*{}**@{},
 <-1.8mm,2.8mm>*{};<0mm,4.9mm>*{}**@{-},
 <-2.8mm,2.9mm>*{};<-4.6mm,4.9mm>*{}**@{-},
   <0.49mm,0.49mm>*{};<2.7mm,2.3mm>*{^2}**@{},
   <-1.8mm,2.8mm>*{};<0.4mm,5.3mm>*{^1}**@{},
   <-2.8mm,2.9mm>*{};<-5.1mm,5.3mm>*{^3}**@{},
 \end{xy}
\ + \
\begin{xy}
 <0mm,0mm>*{\circ};<0mm,0mm>*{}**@{},
 <0mm,-0.49mm>*{};<0mm,-3.0mm>*{}**@{-},
 <0.49mm,0.49mm>*{};<1.9mm,1.9mm>*{}**@{-},
 <-0.5mm,0.5mm>*{};<-1.9mm,1.9mm>*{}**@{-},
 <-2.3mm,2.3mm>*{\circ};<-2.3mm,2.3mm>*{}**@{},
 <-1.8mm,2.8mm>*{};<0mm,4.9mm>*{}**@{-},
 <-2.8mm,2.9mm>*{};<-4.6mm,4.9mm>*{}**@{-},
   <0.49mm,0.49mm>*{};<2.7mm,2.3mm>*{^1}**@{},
   <-1.8mm,2.8mm>*{};<0.4mm,5.3mm>*{^3}**@{},
   <-2.8mm,2.9mm>*{};<-5.1mm,5.3mm>*{^2}**@{},
 \end{xy}\ =\ 0,
 \ \ \ \
\Ba{c}\begin{xy}
 <0mm,0mm>*{\bu};<0mm,0mm>*{}**@{},
 <0mm,0.69mm>*{};<0mm,3.0mm>*{}**@{-},
 <0.39mm,-0.39mm>*{};<2.4mm,-2.4mm>*{}**@{-},
 <-0.35mm,-0.35mm>*{};<-1.9mm,-1.9mm>*{}**@{-},
 <-2.4mm,-2.4mm>*{\bu};<-2.4mm,-2.4mm>*{}**@{},
 <-2.0mm,-2.8mm>*{};<0mm,-4.9mm>*{}**@{-},
 <-2.8mm,-2.9mm>*{};<-4.7mm,-4.9mm>*{}**@{-},
    <0.39mm,-0.39mm>*{};<3.3mm,-4.0mm>*{^3}**@{},
    <-2.0mm,-2.8mm>*{};<0.5mm,-6.7mm>*{^2}**@{},
    <-2.8mm,-2.9mm>*{};<-5.2mm,-6.7mm>*{^1}**@{},
 \end{xy}\Ea
\ + \
 \Ba{c}\begin{xy}
 <0mm,0mm>*{\bu};<0mm,0mm>*{}**@{},
 <0mm,0.69mm>*{};<0mm,3.0mm>*{}**@{-},
 <0.39mm,-0.39mm>*{};<2.4mm,-2.4mm>*{}**@{-},
 <-0.35mm,-0.35mm>*{};<-1.9mm,-1.9mm>*{}**@{-},
 <2.4mm,-2.4mm>*{\bu};<-2.4mm,-2.4mm>*{}**@{},
 <2.0mm,-2.8mm>*{};<0mm,-4.9mm>*{}**@{-},
 <2.8mm,-2.9mm>*{};<4.7mm,-4.9mm>*{}**@{-},
    <0.39mm,-0.39mm>*{};<-3mm,-4.0mm>*{^1}**@{},
    <-2.0mm,-2.8mm>*{};<0mm,-6.7mm>*{^2}**@{},
    <-2.8mm,-2.9mm>*{};<5.2mm,-6.7mm>*{^3}**@{},
 \end{xy}\Ea
 =0,
\ \ \ \ \ \
 \begin{xy}
 <0mm,2.47mm>*{};<0mm,0.12mm>*{}**@{-},
 <0.5mm,3.5mm>*{};<2.2mm,5.2mm>*{}**@{-},
 <-0.48mm,3.48mm>*{};<-2.2mm,5.2mm>*{}**@{-},
 <0mm,3mm>*{\circ};<0mm,3mm>*{}**@{},
  <0mm,-0.8mm>*{\bu};<0mm,-0.8mm>*{}**@{},
<-0.39mm,-1.2mm>*{};<-2.2mm,-3.5mm>*{}**@{-},
 <0.39mm,-1.2mm>*{};<2.2mm,-3.5mm>*{}**@{-},
     <0.5mm,3.5mm>*{};<2.8mm,5.7mm>*{^2}**@{},
     <-0.48mm,3.48mm>*{};<-2.8mm,5.7mm>*{^1}**@{},
   <0mm,-0.8mm>*{};<-2.7mm,-5.2mm>*{^1}**@{},
   <0mm,-0.8mm>*{};<2.7mm,-5.2mm>*{^2}**@{},
\end{xy}=0
$$
A differential in $\cA$ is defined by
$$
d\, \begin{xy}
 <0mm,0.66mm>*{};<0mm,3mm>*{}**@{-},
 <0.39mm,-0.39mm>*{};<2.2mm,-2.2mm>*{}**@{-},
 <-0.35mm,-0.35mm>*{};<-2.2mm,-2.2mm>*{}**@{-},
 <0mm,0mm>*{\bu};<0mm,0mm>*{}**@{},
\end{xy}=0, \ \ \ \
d\, \begin{xy}
 <0mm,-0.55mm>*{};<0mm,-2.5mm>*{}**@{-},
 <0.5mm,0.5mm>*{};<2.2mm,2.2mm>*{}**@{-},
 <-0.48mm,0.48mm>*{};<-2.2mm,2.2mm>*{}**@{-},
 <0mm,0mm>*{\circ};<0mm,0mm>*{}**@{},
 \end{xy} =0, \ \ \ \
d\, \begin{xy}
 <0mm,-0.55mm>*{};<0mm,-3mm>*{}**@{-},
 <0mm,0.5mm>*{};<0mm,3mm>*{}**@{-},
 <0mm,0mm>*{\bu};<0mm,0mm>*{}**@{},
 \end{xy}=
 \Ba{c}
\begin{xy}
 <0mm,-1.3mm>*{};<0mm,-3.5mm>*{}**@{-},
 <0.4mm,0.0mm>*{};<2.4mm,2.1mm>*{}**@{-},
 <-0.38mm,-0.2mm>*{};<-2.8mm,2.5mm>*{}**@{-},
<0mm,-0.8mm>*{\circ};<0mm,0.8mm>*{}**@{},
 <2.96mm,2.4mm>*{\circ};<2.45mm,2.35mm>*{}**@{},
 <2.4mm,2.8mm>*{};<0mm,5mm>*{}**@{-},
 <3.35mm,2.9mm>*{};<5.5mm,5mm>*{}**@{-},
<-2.8mm,2.5mm>*{};<0mm,5mm>*{\ast}**@{},
<-2.8mm,2.5mm>*{};<0mm,5mm>*{}**@{-},
<2.96mm,5mm>*{};<2.96mm,7.5mm>*{\ast}**@{},
<0.2mm,5.1mm>*{};<2.8mm,7.5mm>*{}**@{-},
<0.2mm,5.1mm>*{};<2.8mm,7.5mm>*{}**@{-},
<5.5mm,5mm>*{};<2.8mm,7.5mm>*{}**@{-},
<2.9mm,7.5mm>*{};<2.9mm,10.5mm>*{}**@{-},
    \end{xy}
    \vspace{4mm}
\Ea
\ \ \mbox{where}\ \ \
\Ba{c}\xy
<0mm,0.66mm>*{};<0mm,3mm>*{}**@{-},
 <0.39mm,-0.39mm>*{};<2.2mm,-2.2mm>*{}**@{-},
 <-0.35mm,-0.35mm>*{};<-2.2mm,-2.2mm>*{}**@{-},
 <0mm,0mm>*{\ast};<0mm,0mm>*{}**@{},
   <0.39mm,-0.39mm>*{};<2.9mm,-4mm>*{^2}**@{},
   <-0.35mm,-0.35mm>*{};<-2.8mm,-4mm>*{^1}**@{},
    \endxy
    \Ea:=
    \Ba{c}\xy
<0mm,0.66mm>*{};<0mm,3mm>*{}**@{-},
 <0.39mm,-0.39mm>*{};<2.2mm,-2.2mm>*{}**@{-},
 <-0.35mm,-0.35mm>*{};<-2.2mm,-2.2mm>*{}**@{-},
 <0mm,0mm>*{\bu};<0mm,0mm>*{}**@{},
   <0.39mm,-0.39mm>*{};<2.9mm,-4mm>*{^2}**@{},
   <-0.35mm,-0.35mm>*{};<-2.8mm,-4mm>*{^1}**@{},
    \endxy
    \Ea
    +
   \Ba{c}\xy
<0mm,0.66mm>*{};<0mm,3mm>*{}**@{-},
 <0.39mm,-0.39mm>*{};<2.2mm,-2.2mm>*{}**@{-},
 <-0.35mm,-0.35mm>*{};<-2.2mm,-2.2mm>*{}**@{-},
 <0mm,0mm>*{\bu};<0mm,0mm>*{}**@{},
   <0.39mm,-0.39mm>*{};<2.9mm,-4mm>*{^1}**@{},
   <-0.35mm,-0.35mm>*{};<-2.8mm,-4mm>*{^2}**@{},
    \endxy
    \Ea
$$
Note that the generator $\Ba{c}\xy
<0mm,0.66mm>*{};<0mm,3mm>*{}**@{-},
 <0.39mm,-0.39mm>*{};<2.2mm,-2.2mm>*{}**@{-},
 <-0.35mm,-0.35mm>*{};<-2.2mm,-2.2mm>*{}**@{-},
 <0mm,0mm>*{\ast};<0mm,0mm>*{}**@{},
   <0.39mm,-0.39mm>*{};<2.9mm,-4mm>*{^2}**@{},
   <-0.35mm,-0.35mm>*{};<-2.8mm,-4mm>*{^1}**@{},
    \endxy
    \Ea$ satisfies the second relation in  (\ref{R for LieB}) so that we have a canonical injection of dg props
$$
i: \cU(\gr\cP) \lon \cA
$$
It is easy to see that image of  $\cU(\gr\cP)$ under this injection is a {\em direct}\, summand
 in the complex $(\cA, d)$. Hence {\sc Claim 2} is proven once we show that the cohomology
 of the prop $\cA$ is concentrated in $s$-degree zero.
 \sip


 Using the associativity relation for the generator
 $ \begin{xy}
 <0mm,0.66mm>*{};<0mm,3mm>*{}**@{-},
 <0.39mm,-0.39mm>*{};<2.2mm,-2.2mm>*{}**@{-},
 <-0.35mm,-0.35mm>*{};<-2.2mm,-2.2mm>*{}**@{-},
 <0mm,0mm>*{\bu};<0mm,0mm>*{}**@{},
\end{xy}$ and the Jacobi relation for the generator $\Ba{c}
 \begin{xy}
 <0mm,-0.55mm>*{};<0mm,-2.5mm>*{}**@{-},
 <0.5mm,0.5mm>*{};<2.2mm,2.2mm>*{}**@{-},
 <-0.48mm,0.48mm>*{};<-2.2mm,2.2mm>*{}**@{-},
 <0mm,0mm>*{\circ};<0mm,0mm>*{}**@{},
 \end{xy}\Ea$ one obtains an equality

\Beqrn
 \Ba{c}
\begin{xy}
 <0mm,-1.3mm>*{};<0mm,-3.5mm>*{}**@{-},
 <0.4mm,0.0mm>*{};<2.4mm,2.1mm>*{}**@{-},
 <-0.38mm,-0.2mm>*{};<-2.8mm,2.5mm>*{}**@{-},
<0mm,-0.8mm>*{\circ};<0mm,0.8mm>*{}**@{},
 <2.96mm,2.4mm>*{\circ};<2.45mm,2.35mm>*{}**@{},
 <2.4mm,2.8mm>*{};<0mm,5mm>*{}**@{-},
 <3.35mm,2.9mm>*{};<5.5mm,5mm>*{}**@{-},
<-2.8mm,2.5mm>*{};<0mm,5mm>*{\ast}**@{},
<-2.8mm,2.5mm>*{};<0mm,5mm>*{}**@{-},
<2.96mm,5mm>*{};<2.96mm,7.5mm>*{\ast}**@{},
<0.2mm,5.1mm>*{};<2.8mm,7.5mm>*{}**@{-},
<0.2mm,5.1mm>*{};<2.8mm,7.5mm>*{}**@{-},
<5.5mm,5mm>*{};<2.8mm,7.5mm>*{}**@{-},
<2.9mm,7.5mm>*{};<2.9mm,10.5mm>*{}**@{-},
    \end{xy}
    \vspace{4mm}
    \Ea
    &=&
\Ba{c}\begin{xy}
 <0mm,-1.3mm>*{};<0mm,-3.5mm>*{}**@{-},
 <0.4mm,0.0mm>*{};<2.4mm,2.1mm>*{}**@{-},
 <-0.38mm,-0.2mm>*{};<-2.8mm,2.5mm>*{}**@{-},
<0mm,-0.8mm>*{\circ};<0mm,0.8mm>*{}**@{},
 <2.96mm,2.4mm>*{\circ};<2.45mm,2.35mm>*{}**@{},
 <2.4mm,2.8mm>*{};<0mm,5mm>*{}**@{-},
 <3.35mm,2.9mm>*{};<5.5mm,5mm>*{}**@{-},
<-2.8mm,2.5mm>*{};<0mm,5mm>*{\ast}**@{},
<-2.8mm,2.5mm>*{};<0mm,5mm>*{}**@{-},
<2.96mm,5mm>*{};<2.96mm,7.5mm>*{\bu}**@{},
<0.2mm,5.1mm>*{};<2.8mm,7.5mm>*{}**@{-},
<0.2mm,5.1mm>*{};<2.8mm,7.5mm>*{}**@{-},
<5.5mm,5mm>*{};<2.8mm,7.5mm>*{}**@{-},
<2.9mm,7.5mm>*{};<2.9mm,10.5mm>*{}**@{-},
    \end{xy}
    \vspace{4mm}
\Ea
\ + \
\Ba{c}\begin{xy}
 <0mm,-1.3mm>*{};<0mm,-3.5mm>*{}**@{-},
 <0.4mm,0.0mm>*{};<2.4mm,2.1mm>*{}**@{-},
 <-0.38mm,-0.2mm>*{};<-2.8mm,2.5mm>*{}**@{-},
<0mm,-0.8mm>*{\circ};<0mm,0.8mm>*{}**@{},
 <2.96mm,2.4mm>*{\circ};<2.45mm,2.35mm>*{}**@{},
 <2.4mm,2.8mm>*{};<0mm,5mm>*{}**@{-},
 <3.35mm,2.9mm>*{};<5.5mm,5mm>*{}**@{-},
<-2.8mm,2.5mm>*{};<0mm,5mm>*{\bu}**@{},
<-2.8mm,2.5mm>*{};<0mm,5mm>*{}**@{-},
<2.96mm,5mm>*{};<2.96mm,7.5mm>*{\bu}**@{},
<0.2mm,5.1mm>*{};<2.8mm,7.5mm>*{}**@{-},
<0.2mm,5.1mm>*{};<2.8mm,7.5mm>*{}**@{-},
<5.5mm,5mm>*{};<2.8mm,7.5mm>*{}**@{-},
<2.9mm,7.5mm>*{};<2.9mm,10.5mm>*{}**@{-},
    \end{xy}
    \vspace{4mm}
\Ea
\ - \ 2
\Ba{c}\begin{xy}
 <0mm,-1.3mm>*{};<0mm,-3.5mm>*{}**@{-},
 <-0.4mm,0.0mm>*{};<-2.4mm,2.1mm>*{}**@{-},
 <0.38mm,-0.2mm>*{};<5mm,5mm>*{}**@{-},
 <5mm,5mm>*{};<0mm,10mm>*{}**@{-},
<0mm,-0.8mm>*{\circ};
 <-2.96mm,2.4mm>*{\circ};
 <-2.4mm,2.8mm>*{};<-0.2mm,5mm>*{}**@{-},
 <-3.35mm,2.9mm>*{};<-5.5mm,5mm>*{}**@{-},
<0mm,10.5mm>*{\bu}**@{},
<-2.8mm,7.5mm>*{};<0mm,10.5mm>*{}**@{-},
<-2.96mm,7.5mm>*{\bu}**@{},
<-0.2mm,5.0mm>*{};<-2.8mm,7.5mm>*{}**@{-},
<-0.2mm,5.1mm>*{};<-2.8mm,7.5mm>*{}**@{-},
<-5.5mm,5mm>*{};<-2.8mm,7.5mm>*{}**@{-},
<0mm,10.5mm>*{};<0mm,13.5mm>*{}**@{-},
    \end{xy}
    \vspace{4mm}
\Ea
\\
&=&
-3\
\Ba{c}\begin{xy}
 <0mm,-1.3mm>*{};<0mm,-3.5mm>*{}**@{-},
 <-0.4mm,0.0mm>*{};<-2.4mm,2.1mm>*{}**@{-},
 <0.38mm,-0.2mm>*{};<2.5mm,2.5mm>*{}**@{-},
 <2.5mm,7.5mm>*{};<0mm,10mm>*{}**@{-},
<0mm,-0.8mm>*{\circ};
 <-2.96mm,2.4mm>*{\circ};
 <-2.4mm,2.8mm>*{};<2.5mm,7.5mm>*{}**@{-},
 <-3.35mm,2.9mm>*{};<-5.5mm,5mm>*{}**@{-},
<0mm,10.5mm>*{\bu}**@{},
<-2.8mm,7.5mm>*{};<0mm,10.5mm>*{}**@{-},
<-2.96mm,7.5mm>*{\bu}**@{},
<2.5mm,2.5mm>*{};<-2.8mm,7.5mm>*{}**@{-},
<-5.5mm,5mm>*{};<-2.8mm,7.5mm>*{}**@{-},
<0mm,10.5mm>*{};<0mm,13.5mm>*{}**@{-},
    \end{xy}
    \vspace{4mm}
    \Ea= -3\ \frac{\Ba{c}\xy <0mm,0mm>*{\bu};<0mm,0mm>*{}**@{},
 <0mm,0.69mm>*{};<0mm,3.0mm>*{}**@{-},
 <0.39mm,-0.39mm>*{};<2.4mm,-2.4mm>*{}**@{-},
 <-0.35mm,-0.35mm>*{};<-1.9mm,-1.9mm>*{}**@{-},
 <-2.4mm,-2.4mm>*{\bu};<-2.4mm,-2.4mm>*{}**@{},
 <-2.0mm,-2.8mm>*{};<0mm,-4.9mm>*{}**@{-},
 <-2.8mm,-2.9mm>*{};<-4.7mm,-4.9mm>*{}**@{-},
    <0.39mm,-0.39mm>*{};<3.3mm,-4.0mm>*{^3}**@{},
    <-2.0mm,-2.8mm>*{};<0.5mm,-6.7mm>*{^2}**@{},
    <-2.8mm,-2.9mm>*{};<-5.2mm,-6.7mm>*{^1}**@{},
 \endxy\Ea}{\Ba{c}\begin{xy}
 <0mm,0mm>*{\circ};<0mm,0mm>*{}**@{},
 <0mm,-0.49mm>*{};<0mm,-3.0mm>*{}**@{-},
 <0.49mm,0.49mm>*{};<1.9mm,1.9mm>*{}**@{-},
 <-0.5mm,0.5mm>*{};<-1.9mm,1.9mm>*{}**@{-},
 <-2.3mm,2.3mm>*{\circ};<-2.3mm,2.3mm>*{}**@{},
 <-1.8mm,2.8mm>*{};<0mm,4.9mm>*{}**@{-},
 <-2.8mm,2.9mm>*{};<-4.6mm,4.9mm>*{}**@{-},
   <0.49mm,0.49mm>*{};<2.7mm,2.3mm>*{^2}**@{},
   <-1.8mm,2.8mm>*{};<0.4mm,5.3mm>*{^3}**@{},
   <-2.8mm,2.9mm>*{};<-5.1mm,5.3mm>*{^1}**@{},
 \end{xy}\Ea}
\Eeqrn
where the horizontal line stands for the properadic composition in accordance with the labels shown.

\sip
This result propmts us to consider a dg associative non-commutative algebra $\mathbf{A}_n$ generated by
degree zero variable $\{x_i\}_{1\leq i\leq n}$ and degree $-1$ generators $\{u_{i,i+1,u+2}\}_{1\leq i\leq n-2}$
with the differential
$$
d u_{i,i+1,i+2}= -3[[x_i,x_{i+2}],x_{i+1}]
$$
or, equivalently (after rescaling the generators $u_\bu$),
with the differential
$$
d u_{i,i+1,i+2}= [[x_i,x_{i+2}],x_{i+1}]
$$
Arguing in exactly the same way as in \cite{CMW} one concludes that the cohomology of the dg operad
operad $\cA$ is concentrated in $s$-degree zero if and only if the collections of algebras
$\mathbf{A}_n$, $n\geq 3$, has cohomology concentrated in ordinary degree zero. The latter fact is established
in the appendix.
The proof is completed.
\end{proof}

\subsection{Representations of  $\HoLonB$ and quantizable Poisson structures} Let $V=\R^d$ be a $d$-dimensional vector space, and $\f_V=\prod_{n\geq 1} \odot^n V^*$  the commutative algebra of formal power series functions on $V$, with the obvious complete topology. If $\mathrm{Der}(\f_V)$ stands for the Lie algebra of continuous derivations
of $\f_V$, then the Lie algebra of formal polyvector fields on $V$ is defined as the Lie algebra of continuous multiderivations,
$$
\cT_{poly}(V):= \widehat{\odot^\bu_{\f_V}}\left(\mathrm{Der}(\f_V)[1]\right)\cong \prod_{m\geq 0} \wedge^m V\ot \f_{V} \simeq \prod_{m,n\geq 0} \wedge^m V \ot \odot^n V^* \, .
$$
There is an obvious chain of injections of topological commutative algebras,
 $$
 \ldots \lon \f_{\R^d} \lon  \f_{\R^{d+1}} \lon  \f_{\R^{d+2}} \lon  \ldots.
 $$
 We denote the associated {\em direct}\, limit by
 $$
 \f_{\R^\infty} := \lim_{d\rar \infty}  \f_{\R^d}.
 $$
For $V=\R^\infty$ we define $\cT_{poly}(V)$ as the Lie algebra of continuous multiderivations of $\f_{\R^\infty}$, i.e.,
$$
\cT_{poly}(V)\cong \prod_{m\geq 0} \Hom(\wedge^m \R^\infty, \f_{\R^\infty}) \, .
$$
We can also consider the space  $\cT_{poly}(V)[[\hbar]]$ of formal power series in a formal variable $\hbar$ with coefficients in $\cT_{poly}(V)$. These arguments can be easily generalized to a finite/infinite
dimensional {\em graded}\, vector space $V$.

\sip

Consider now a representation of our minimal resolution
$$
\rho: \HoLonB \lon \cE nd_V
$$
in a (possible, infinite-dimensional) dg vector space $V$. It is uniquely determined
by the values of $\rho$ on the generators of $\HoLonB$,
$$
\rho\left(
\resizebox{15mm}{!}{
\xy
(-9,-6)*{};
(0,0)*+{a}*\cir{}
**\dir{-};
(-5,-6)*{};
(0,0)*+{a}*\cir{}
**\dir{-};
(9,-6)*{};
(0,0)*+{a}*\cir{}
**\dir{-};
(5,-6)*{};
(0,0)*+{a}*\cir{}
**\dir{-};
(0,-6)*{\ldots};
(-10,-8)*{_1};
(-6,-8)*{_2};
(10,-8)*{_n};
(-9,6)*{};
(0,0)*+{a}*\cir{}
**\dir{-};
(-5,6)*{};
(0,0)*+{a}*\cir{}
**\dir{-};
(9,6)*{};
(0,0)*+{a}*\cir{}
**\dir{-};
(5,6)*{};
(0,0)*+{a}*\cir{}
**\dir{-};
(0,6)*{\ldots};
(-10,8)*{_1};
(-6,8)*{_2};
(10,8)*{_m};
\endxy}\right):=\pi^m_n(a)\in \wedge^m V \ot \odot^n V^*.
$$
We can assemble these values into a formal power series
$$
\pi^\diamond := \sum_{m,n\geq 0}\sum_{a\geq 0}\hbar^a \pi^m_n(a)\in \cT_{poly}(V)
$$
which gives us a formal polyvector field on $V$. The values $\pi^m_n(a)$ can not be chosen arbitrarily as the map $\rho$ must respect differentials in $\HoLonB$ and $V$,
$$
\rho\circ \delta_\diamond= d\circ \rho.
$$
Untwisting the definition of $\delta_\diamond$, we conclude that the above formal power series
$\pi^\diamond$ (with  $\pi^1_1(0):=d$) comes from a representation of $\HoLonB$ if and only if it satisfies the equation
$$
\frac{1}{2}[\pi^\diamond,\pi^\diamond]_2 + \frac{\hbar}{4!} [\pi^\diamond,\pi^\diamond,\pi^\diamond,\pi^\diamond]_4
+ \ldots =0,
$$
where the collection of operators,
$$
\left\{[\, \ , \ldots ,\ ]_{2n}: \cT_{poly}(V)^{\ot 2n}\rar \cT_{poly}(V)[3-4n] \right\}_{n\geq 1}
$$
comes from the values on the graphs $\Upsilon_{2n}$ from \S 3.1 of the standard morphism \cite{Wi1}
of dg Lie algebras
$$
\fGCor_2 \rar CE^\bu( \cT_{poly}(V),  \cT_{poly}(V)),
$$
$CE^\bu( \cT_{poly}(V),  \cT_{poly}(V))$ being the Chevalley-Eilenberg deformation complex of the Lie algebra of polyvector fields. Therefore formal {\em quantizable}\, Poisson structures on a graded vector space $V$ (viewed as a formal manifold) come from representations of our properad $\HoLonB$.
There are plenty of examples of such quantizable Poisson structures on {\em finite}-dimensional vector spaces, one for each ordinary formal graded Poisson structure $\pi$ on $V$ (which is, by definition, an element of $ \cT_{poly}(V)$ which satisfies the standard Schouten equation $[\pi,\pi]_2=0$). However
the association
$$
\pi \lon \pi^\diamond
$$
is highly non-trivial and depends on the choice of a Drinfeld associator \cite{MW2}. It is an open problem to find a non-trivial example of a quantizable Poisson structure in {\em infinite}\, dimensions. Perhaps, for any graded vector space $V$ equipped with an odd symplectic form, the associated total space of cyclic words
$$
V:= \Pi_{n\geq 1} (\ot^n W)_{\Z_n}
$$
comes equipped with such a structure given by formulae from Theorem 4.3.3 in \cite{MW0}; however it is hard to check this conjecture by a direct computation as it involves infinitely many equations.

\bip

\bip


\appendix

\section{}

\subsection{\bf Lemma}
{}
{\label{lem::free::grob}}
{\em
Let $\{c_\sigma| \sigma\in S_3\}$ be a collection of $6$ numbers such that
\begin{equation}
\begin{split}
\label{eq::prop::lm}
 \text{for each pair $(c_{ijk},c_{ikj})$ with the same first index } i \\
 \text{ at least one of these elements is different from zero
}
\end{split}
\end{equation}
Then the associative algebra
$
A:=
\Bbbk \left\langle
x_1,\ldots, x_n \left|
\begin{array}{c}
\sum_{\sigma\in S_3} c_\sigma x_{i+\sigma(1)} x_{i+\sigma(2)} x_{i+\sigma(3)}, \\
 i = 0,\ldots, n-2
\end{array}
\right.
\right\rangle
$
has global dimension $2$.
}
\begin{proof}
Let us consider any linear ordering of the set of generators, such that
$$
\forall k,l,m \ \ x_{3k} > x_{3l+2}, x_{3m+1}
$$
We extend this ordering to a degree-lexicographical ordering of the set of monomials in the free associative algebra $\Bbbk\langle x_1,\ldots,x_n\rangle$.
The leading monomials of relation $\sum_{\sigma\in S_3} c_\sigma x_{i+\sigma(1)} x_{i+\sigma(2)} x_{i+\sigma(3)}$ are different for all $i$ because they contain different letters $\{x_{i+1},x_{i+2},x_{i+3}\}$.
Moreover, there is exactly one number divisible by $3$ in each subsequent triple of integer numbers, thus after reordering we have
$\{i+1,i+2,i+3\}= \{3s,r,t\}$ for appropriate $r,s$ and $t$, such that $r$ and $t$ are not divisible by $3$.
Recall that by property~\eqref{eq::prop::lm}
at least one of the two monomials $c_{{3s} r t}x_{3s}x_r x_t$ and
$c_{{3s} t r}x_{3s}x_t x_r$ is different from zero. Hence, the first letters in the leading monomials
of the relation $\sum_{\sigma\in S_3} c_\sigma x_{i+\sigma(1)} x_{i+\sigma(2)} x_{i+\sigma(3)}$
have index divisible by $3$ and two remaining letters is not divisible by $3$.
  Consequently, the leading monomials of generating relations have no compositions and the set of generating relations form a \emph{strongly free} Gr\"obner bases following that the algebra $A$ has global dimension $2$. (See~\cite{ufn} \S4.3 for details on strongly free relations.)
\end{proof}

\subsection{\bf Corollary}
{\em
The minimal free resolution $\mathbf{ A}_{n}$ of the algebra
$${A}_{n}:= \Bbbk\left\langle x_1,\ldots,x_n \left|
\begin{array}{c}
 [[x_i,x_{i+2}],x_{i+1}] \\
i=1,\ldots,n-2
\end{array}
\right.
\right\rangle$$
is generated by $x_1,\ldots,x_n$ and  $u_{1,2,3}\ldots,u_{n-2,n-1,n}$ such that
}
$$
\begin{array}{ccc} \deg(x_i)=0, & & \deg(u_{i,i+1,i+2})=-1;\\
   d(x_i)=0; & \quad &
     d(u_{i,i+1,i+2})= [[x_i,x_{i+2}],x_{i+1}].
\end{array}
$$
\begin{proof}
Let us expand the commutators in the relations we are working with:
$$
 [[x_1,x_{3}],x_{2}]
=
(x_1 x_3 x_2 - x_3 x_1 x_2 - x_2 x_1 x_3 + x_2 x_3 x_1)
$$
As we can see they satisfy the condition~\eqref{eq::prop::lm} of Lemma~{\ref{lem::free::grob}} and algebra $A_n$ has global dimension $2$, meaning that the following complex
$$
0\to \mathrm{span}\langle\text{relations}\rangle \otimes A_n \to \mathrm{span}\langle\text{generators}\rangle\otimes A_n \to A_n \to \Bbbk \to 0
$$
is acyclic in the leftmost term and, consequently, acyclic everywhere. Therefore, the minimal resolution of $A_n$ is generated by generators and generating relations of $A_n$.
\end{proof}

\def\cprime{$'$}

\end{document}